\documentclass[conference]{IEEEtran}
\IEEEoverridecommandlockouts
\usepackage{draftwatermark}
\usepackage{cite}
\usepackage{algorithmic}
\usepackage[utf8]{inputenc}
\usepackage{epsfig}
\usepackage{subfigure}
\usepackage{amsfonts,amsthm,amsmath,amssymb,eso-pic}
\usepackage{textcomp}
\usepackage{comment}
\usepackage{gensymb}
\usepackage[font=small,skip=0pt]{caption}
\usepackage[inline]{enumitem}
\usepackage[colorlinks=true,       
            linkcolor=blue,        
            citecolor=blue,        
            urlcolor=blue,         
            bookmarks=true,        
            breaklinks=true        
]{hyperref}

\usepackage[capitalise,nameinlink]{cleveref} 
\usepackage{xcolor}

\hypersetup{
colorlinks=true, 
breaklinks=true, 
urlcolor= black, 
citecolor=black,
linkcolor= blue, 
bookmarksopen=true,
pdftitle={Optimized adaptive MPC for lateral control of autonomous vehicles}, 
pdfauthor={Y. Kebbati},
}

\graphicspath{{Figures/}}
\def\BibTeX{{\rm B\kern-.05em{\sc i\kern-.025em b}\kern-.08em
    T\kern-.1667em\lower.7ex\hbox{E}\kern-.125emX}}
\newcommand{\RNum}[1]{\uppercase\expandafter{\romannumeral #1\relax}}
\newcommand{\linebreakand}{%
  \end{@IEEEauthorhalign}
  \hfill\mbox{}\par
  \mbox{}\hfill\begin{@IEEEauthorhalign}
}
\title{Optimized adaptive MPC for lateral control of autonomous vehicles}

\author{
\IEEEauthorblockN{Yassine Kebbati\textsuperscript{1}, Vicenç Puig\textsuperscript{2}, Naima Ait-Oufroukh\textsuperscript{1}}
\IEEEauthorblockA{\textsuperscript{1}\textit{IBISC-EA4526, Université Paris-Saclay, Evry, France} \\
\textsuperscript{2}\textit{CS2AC, UPC, Barcelona, Spain} \\
\href{mailto:yassine.kebbati@univ-evry.fr}{yassine.kebbati@univ-evry.fr}, 
\href{mailto:vicenc.puig@upc.edu}{vicenc.puig@upc.edu}, 
\href{mailto:naima.aitoufroukh@univ-evry.fr}{naima.aitoufroukh@univ-evry.fr}} 
\linebreak
\IEEEauthorblockN{Vincent Vigneron\textsuperscript{1}, Dalil Ichalal\textsuperscript{1}}
\IEEEauthorblockA{\textsuperscript{1}\textit{IBISC-EA4526, Université Paris-Saclay, Evry, France} \\
\href{mailto:vincent.vigneron@univ-evry.fr}{vincent.vigneron@univ-evry.fr}, 
\href{mailto:dalil.ichalal@univ-evry.fr}{dalil.ichalal@univ-evry.fr}}
}

\begin{document}

\maketitle

\SetWatermarkText{Preprint}
\begin{abstract}

Autonomous vehicles are the upcoming solution to most transportation problems such as safety, comfort and efficiency. The steering control is one of the main important tasks in achieving autonomous driving. Model predictive control (MPC) is among the fittest controllers for this task due to its optimal performance and ability to handle constraints. This paper proposes an adaptive MPC controller (AMPC) for the path tracking task, and an improved PSO algorithm for optimising the AMPC parameters. Parameter adaption is realised online using a lookup table approach. The propose AMPC performance is assessed and compared with the classic MPC and the Pure Pursuit controller through simulations.

\end{abstract}

\section{Introduction}
The increase in road accidents and traffic jams all around the world has encouraged research to further develop autonomous driving technologies \cite{1}. This field has received great attention and has seen some important advances in the last decades \cite{2,3}, especially urban autonomous driving that is still a big challenge today. Among the most famous competitions in the field was the DARPA grand challenge that took place in 2007 \cite{4}, where several research groups tested their autonomous driving systems in an urban-like environment but only six of them managed to finish the race.

Self-driving vehicles do not need human intervention to navigate their path, they are fully controlled by the control module that commands the different actuators such as the accelerator and the steering wheel. Autonomous vehicles are divided into 5 automation levels; level zero represents the conventional cars with no automation, level one includes simple automatic driving systems like adaptive cruise control and electronic stability control. Level two introduces advanced systems like speed and steering control or emergency braking systems. Level three allows the vehicle to sense the environment through multiple sensors and drive autonomously while requiring human intervention in infeasible situations. In level four, the vehicle is fully autonomous with only occasional human intervention and only for certain driving modes, while in level five the vehicle is fully autonomous in all driving modes and without any human intervention \cite{5}.
Path tracking is a fundamental part of autonomous driving, it ensures that the vehicle follows a predefined trajectory and this is achieved by controlling the vehicle lateral dynamics. Various control strategies for this task are reported in the literature. For example, authors of \cite{6} developed a lateral control system based on the adaptive pure pursuit controller, they used a PI controller to reduce the lateral offset and improve the tracking error caused by the look ahead distance. In \cite{7}, an output-feedback robust controller is introduced to assist in the path tracking task, the controller considers different driver's characteristics described by uncertain parameters such as delay time and achieves good stability and tracking performance through regional pole placement. Han \textit{et al}. \cite{8} designed an adaptive neural network PID controller for path tracking, they used it with a second order vehicle model whose parameters are estimated using the forgetting factor least square estimation. Kebbati \textit{et al.} \cite{26} developed a self-adaptive PID controller for speed regulation using neural networks and then genetic algorithms. Hu \textit{et al}. \cite{9} introduced an integral sliding mode controller for path tracking of an independently four-wheel actuated vehicle, the proposed nonlinear feedback controller considers multi-input multi-output (MIMO) systems with time varying trajectories. However, most of the above mentioned strategies cannot handle constraints that are imposed by the safety and physical limitations of actuators.

Model predictive control (MPC) stands out in this field due to the fact that it systematically handles constraints on control, output and state signals, in addition to the fact that it easily handles MIMO systems. In \cite{10}, Zhang \textit{et al}. designed an MPC controller for path following with low computation load, for which they used Laguerre functions to approximate the control signal, and exponential weight to improve the tracking accuracy. Guo \textit{et al}. \cite{11} developed a path-tracking MPC controller that considers measurable disturbances, the authors considered the varying road conditions and small-angle assumptions as a form of measurable disturbance and solved the control problem using the differential evolution algorithm. In MPC control methods, there are many parameters to tune such as the weighting matrices and prediction and control horizons which is not trivial. This has not been addressed in the above mentioned papers, but in \cite{12} and \cite{13} the authors used fuzzy inference systems (FIS) to tune the weighting matrices of the MPC. Papers \cite{14} and \cite{15} used genetic algorithms (GA) and particle swarm optimization (PSO) to optimize the MPC tunable parameters. Nevertheless, the controllers in these approaches remain non adaptive to varying parameters and external uncertainties. In \cite{16}, Lin \textit{et al}. designed an adaptive MPC controller that adapts to estimated varying road friction and cornering stiffness coefficients and in \cite{17} an adaptive MPC was designed for lane keeping systems where the steering offset is continuously learned from measured data and adapted. 

This paper proposes an adaptive MPC controller for the lateral control task, the contributions of this work are threefold: first, the design of an optimized adaptive MPC controller with Laguerre functions for the lateral control task. Second, the optimization of MPC parameters and cost function weight matrices with a new improved PSO algorithm. Third, the adaptation of the MPC, cost function and Laguerre function parameters online using a lookup table. Section \RNum{2} of this paper presents the vehicle lateral dynamics modeling and MPC control formulation. The controller optimization with improved PSO algorithm is introduced in section \RNum{3}. The results of the proposed controller are evaluated and discussed in section \RNum{4}. Section \RNum{5} gives final conclusions and some suggestions for the continuation of future research.

\section{Vehicle Modeling and Controller Design}\label{sec:Vehicle Modeling}

\subsection{Vehicle Lateral Dynamics Model} \label{sec:Vehicle Model}
\begin{figure}[b]
\centering
\includegraphics[width=0.45\textwidth]{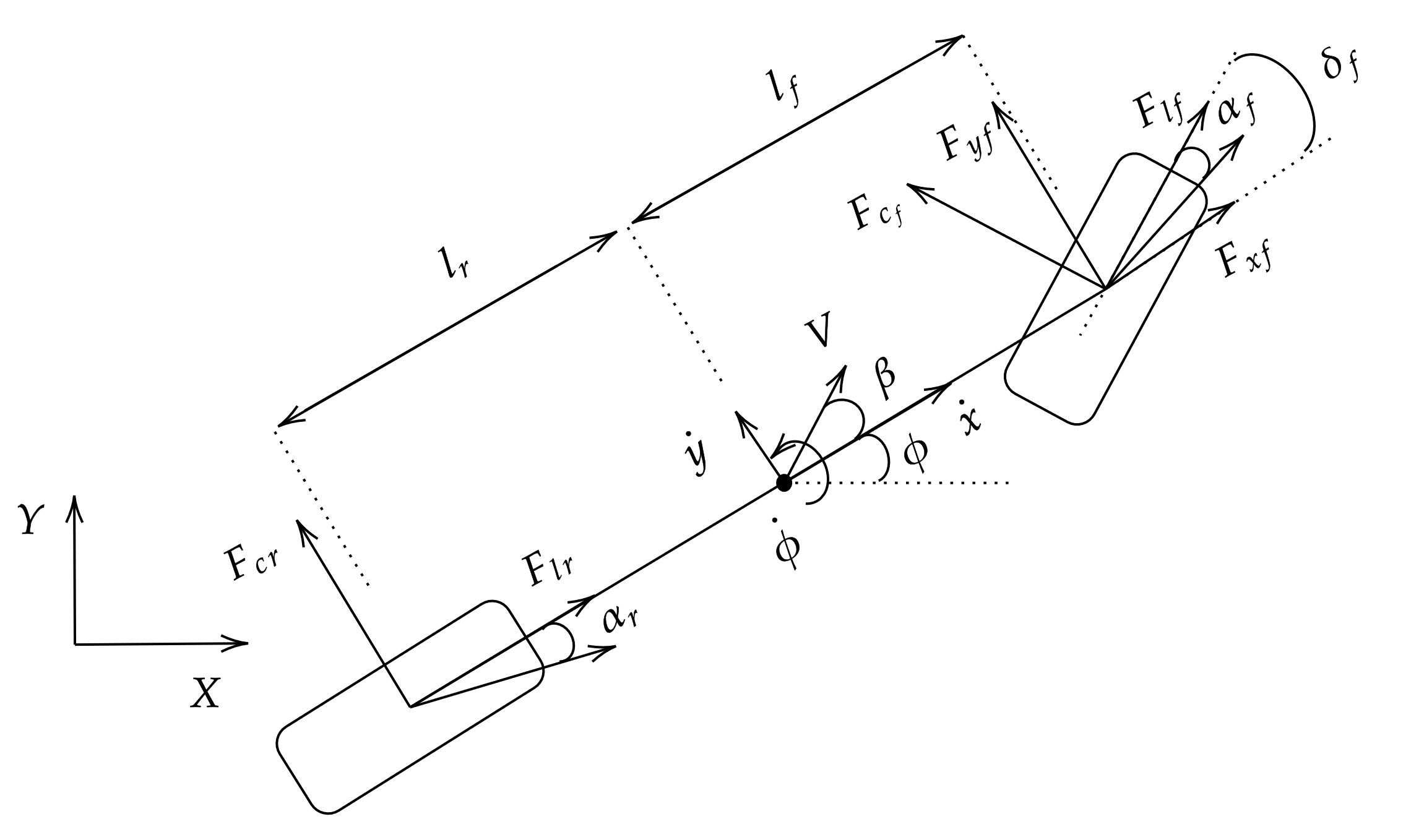}
\caption{2-DOF Bicycle dynamic model.}
\label{fig:1}
\end{figure}

Path tracking is associated with the lateral motion of the vehicle, which is the displacement on the y-axis and the yaw motion around the z-axis. Therefore, the bicycle (2-DOF) dynamics model shown in Fig. \ref{fig:1}, can accurately represent the vehicle motion and it is used to formulate the MPC controller. Applying newton's law of motion along the y-axis and z-axis yields the following equations:
\setlength{\textfloatsep}{2pt}

\begin{equation}
\label{eq1}
\left\{
    \begin{array}{ll}
 m(\Ddot{y}+\Dot{x}\Dot{\phi}) &= 2F_{yf} + 2F_{yr}\\
 I_z \Ddot{\phi} &= 2l_f F_{yf} - 2l_r F_{yr}
    \end{array}
\right.
\end{equation}
where, $m$ represents the mass of the vehicle, $y$, $x$ and $\phi$ are the respective lateral and longitudinal positions and the heading angle of the vehicle, $F_{yf}$ and $F_{yr}$ are the longitudinal forces of the front and rear wheels respectively and $l_f$ and $l_r$ are the distances from the front and rear wheel axles to the vehicle's center of gravity (CG). The longitudinal force $F_y$ is a result of the nonlinear interaction between the tires and the road, In this work the linearized tire model under the assumption of small slip angles is used and it is given as follows:

\begin{equation}
\label{eq2}
\left\{
    \begin{array}{ll}
 F_{yf} &= C_{yf} \alpha_f\\
 F_{yr} &= C_{yr} \alpha_r
    \end{array}
\right.
\end{equation}
with $C_{y(f,r)}$ and $a_{(f,r)}$ being the cornering stiffness coefficient and the lateral slip angle for the front and rear wheels respectively. The lateral slip angles are given by the following equations: 

\begin{equation}
\label{eq3}
\left\{
    \begin{array}{ll}
 \alpha_f &= \delta_f - \gamma_f\\
 \alpha_r &= \delta_r - \gamma_r
    \end{array}
\right.
\end{equation}
where $\delta_{(f,r)}$ is the front and rear wheel steering angles where only the front wheel is considered steerable $(\delta_r = 0)$, $\gamma_{(f,r)}$ define the angles between the direction of the front/rear wheel directions and the longitudinal velocity given as:

\begin{equation}
\label{eq4}
\left\{
    \begin{array}{ll}
 \tan (\gamma_f) &= \frac{\Dot{y} + l_f \Dot{\phi}}{\Dot{x}}\\
 \tan (\gamma_r) &= \frac{\Dot{y} - l_r \Dot{\phi}}{\Dot{x}}
    \end{array}
\right.
\end{equation}
The transformation from the vehicle body frame to the inertial frame is obtained by equation (\ref{eq5}):
\begin{equation}
    \label{eq5}
    \dot{Y} = \dot{x} \sin{\phi} + \dot{y} \cos{\phi}
\end{equation}
Hence, by considering equation (\ref{eq5}), using angle approximation and substituting equations (\ref{eq4}) and (\ref{eq3}) into (\ref{eq2}) and then equation (\ref{eq2}) into (\ref{eq1}) we obtain the lateral dynamics model as the following:

\begin{equation}
\label{eq6}
\left\{
    \begin{array}{ll}
        m(\Ddot{y}+\Dot{x}\Dot{\phi}) &= 2[C_{yf}(\delta_f-\frac{\Dot{y} + l_f \Dot{\phi}}{\Dot{x}}) + C_{yr}\frac{l_r \Dot{\phi}- \Dot{y}}{\Dot{x}}]\\
        I_z \Ddot{\phi} &= 2[l_f C_{yr}(\delta_f - \frac{\Dot{y} + l_f \Dot{\phi}}{\Dot{x}}) - l_r C_{yr} \frac{l_r \Dot{\phi} - \Dot{y}}{\Dot{x}}]
    \end{array}
\right.
\end{equation}
The above mentioned model can be transformed into the following continuous state space representation: 

\begin{equation}
\label{eq7}
\left\{
    \begin{array}{ll}
        \Dot{x}_c &= A_cx_c + B_cu\\
        y_c &= C_cx_c
    \end{array}
\right.
\end{equation}
where $x_c = [\Dot{y}\ \phi\ \Dot{\phi}\ y]^T$ is the state vector, $y_c = y$ is the output of the system, $u = \delta_f$ is the control input and $A_c$, $B_c$ and $C_c$ are the respective state, control and output matrices given by:
$A_c = \left[{\begin{array}{cccc}-2\frac{C_{yf} + C_{yr}}{m \dot{x}} & 0 & -\dot{x} - 2\frac{C_{yf} l_f -C_{yr} l_r}{m \dot{x}} & 0\\ 0 & 0 & 1 & 0\\-2 \frac{C_{yf} l_f - C_{yr} l_r}{I_z \dot{x}} & 0 & -2\frac{C_{yf} l_f^2+ C_{yr} l_r^2}{I_z \dot{x}} & 0\\1 & \dot{x} & 0 & 0\end{array}}\right],$\\
$B_c = \left[\begin{array}{c} 2\frac{C_{yf}}{m} \\ 0 \\ 2\frac{C_{yf} l_f}{I_z}\\ 0 \end{array}\right],$ \hspace{5mm}
$C_c = \left[\begin{array}{cccc} 0 &  0 & 0 & 1\end{array}\right].$

\subsection{MPC Controller Design}\label{sec:mpc}
\begin{figure}[b]
\centering
\includegraphics[width=0.4\textwidth]{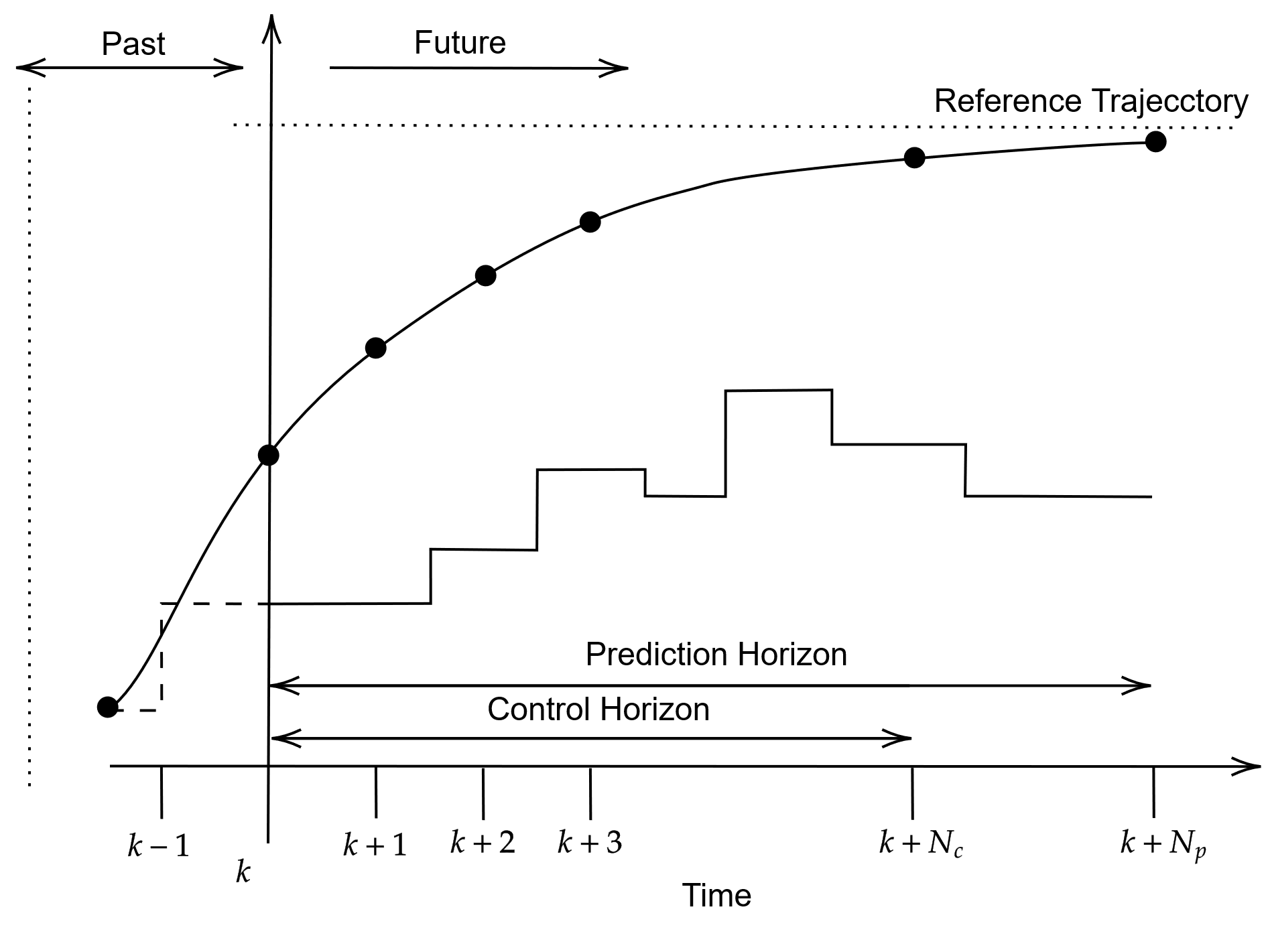}
\caption{MPC receding horizon concept.}
\label{fig:2}
\end{figure}

The linear MPC controller with Laguerre functions is considered in this paper due to its high computational efficiency. The idea of MPC control is to use a descriptive model of the plant to predict its behaviour in the future along a prediction horizon $N_p$, and produce an optimal control sequence along the control horizon $N_c$ that minimizes the error between the set-point and the output of the plant. In the context of path tracking, MPC aims to minimize the gap between the reference trajectory and the predicted trajectory. The optimal control sequence is achieved by minimizing a quadratic cost function subject to constraints, therefore, solving a constrained convex optimization problem. However, only the first term of this sequence is applied at each iteration and the whole prediction-optimization process is repeated iteratively, this is known as the receding horizon principle (Fig. \ref{fig:2}). The linearized vehicle dynamic model of (\ref{eq7}) is discretized to be used as the prediction model of the MPC controller:

\begin{equation}
\label{eq8}
\left\{
    \begin{array}{ll}
        x_{(k+1)} &= A_kx_{(k)} + Bu_{(k)}\\
        y_{(k)} &= Cx_{(k)}
    \end{array}
\right.
\end{equation}

To allow the MPC controller to be adaptive to varying longitudinal speeds, the state matrix $A_k$ is iteratively updated with the current longitudinal velocity to produce an adaptive prediction model for varying speed profiles and not only fixed longitudinal speeds as is the case in most studies in the literature \cite{10,12,16,17}. This standard model is then augmented by embedding an integrator and introducing the output $y_{(k)}$ to the state vector as follows:

\begin{equation}
\label{eq9}
\left\{
    \begin{array}{ll}
        \Tilde{x}_{(k+1)} &= \Tilde{A}_k  \Tilde{x}_{(k)} + \Tilde{B} \Delta \Tilde{u}_{(k)}\\
       \Tilde{y}_{(k)} &= \Tilde{C} \Tilde{x}_{(k)}
    \end{array}
\right.
\end{equation}
where the new state $\Tilde{x} = \left[\Delta x_{(k)}\  y_{(k)}\right]^T$ is connected to the output, the new input signal is $\Delta \Tilde{u}_{(k)}$ and the augmented system matrices are given by:\vspace{0.25cm}

$\Tilde{A}=
\left[\begin{array}{cc} A_k & o_m^T \\CA_k & I_{q \times q}  \end{array}\right]$,
$\Tilde{B} = \left[\begin{array}{c} B \\ CB\end{array}\right]$,
$\Tilde{C}=\left[o_m\ I_{q\times q}\right].$
$o_m =\overbrace{[0\ 0...0]}^\text{$n$}$ is a zero vector and $I_{q \times q}$ is an identity matrix with $n$ being the number of states, $m$ the number of inputs and $q$ the number of outputs. At $k_i>0$ the state vector $x(k)$ is assumed available providing the current plant information, the future control trajectory along the control horizon $N_c$ is given by:
$$\Delta u(k_i), \Delta u(k_i +1),..., \Delta u(k_i + N_c -1)$$
The future plant information is predicted along the prediction horizon $N_p$ using the model given in equation (\ref{eq9}) and the future control parameters through the following equation:

\begin{equation}
\label{eq10}
Y = F x_{(k_i)} + \Phi \Delta U
\end{equation}
where $Y$ is the future outputs, $F$ is the future prediction matrix, $\Phi$ is the future control matrix and $\Delta U$ is the future control sequence given by:
$$Y=[y(k_i+1|k_i)\ y(k_i+2|k_i y(k_i+3|k_i)\hdots y(k_i+N_p|k_i)]^T $$
$$\Delta U=[\Delta u(k_i)\ \Delta u(k_i+1)\ \Delta u(k_i+2) \hdots \Delta u(k_i+N_c-1)]^T $$

$$F = \left[\begin{array}{c} CA \\CA^2\\CA^3\\\vdots \\CA^{N_p}  \end{array}\right]$$ 
$$\Phi = \left[\begin{array}{cccc} CB&0&\hdots&0 \\CAB&CB&\hdots&0\\CA^2B&CAB&\hdots&0\\ \vdots&\vdots&\vdots&\vdots \\CA^{N_p-1}B&CA^{N_p-2}B&\hdots& CA^{N_p-N_c}B \end{array}\right]$$
\setlength{\textfloatsep}{2pt}

The goal of MPC is to find the best $\Delta U$ that minimizes the error $(R_s-Y)$, where $R_s = [1\ 1 \hdots 1]r(k_i)$ is the reference trajectory containing $N_p$ set-point information $r(k_i)$. Then, the constrained optimization problem is formulated as: 
\begin{align}
\label{eq11}
min\ &J = (R_s - Y)^T Q (R_s - Y) + \Delta U^T R \Delta U\\
\label{eq12}
s.t:\ & x(k+1) = A x(k) + B \Delta u(k) \\
\label{eq13}
    & \Delta u_{min}\leq \Delta u \leq \Delta u_{max}\\
\label{eq14}
    & u_{min}\leq u \leq u_{max}\\
\label{eq15}
    & y_{min}\leq y \leq y_{max}
\end{align}
\setlength{\textfloatsep}{4pt}\noindent
Equation (\ref{eq11}) is the cost function to be minimized, (\ref{eq12}) is the vehicle dynamic model and equations (\ref{eq13}) to (\ref{eq15}) represent the constraints on the control, control increment and output signals. $Q$ and $R$ are weighting matrices that penalise the tracking precision and control magnitude respectively.
The set of constraints can be reformulated and expressed in terms of the decision variable $\Delta U$ as:
\setlength{\textfloatsep}{2pt}
\begin{equation}
\label{eq16}
M \Delta U \leq \gamma
\end{equation}
such that $M$ contains 3 reformulation sub-matrices and $\gamma$ contains the set of upper and lower bounds on the control, control increment and output signals. The model predictive control becomes a quadratic programming (QP) problem with quadratic cost function and linear inequality constraints formulated in the following form:
\begin{equation}
\label{eq17}
\left\{
    \begin{array}{ll}
        &J  = \frac{1}{2}x^TEx +x^TK \\
         &Mx \leqslant \gamma
    \end{array}
\right.
\end{equation}
where $E, K, M$ and $\gamma$ are compatible matrices and vectors and $x$ is the decision variable ($\Delta u$) in the quadratic programming problem, $E$ is assumed symmetric and positive definite. \vspace{0.3cm}

\subsubsection{Laguerre function approximation}

The Laguerre function is used to reduce the computation burden by approximating the control sequence $\Delta U$, the orthogonal feature of the Laguerre function allows to realize long control horizons without using a large number of parameters. The discrete form of a Laguerre function is given by:
\begin{equation}\label{18}
\begin{split}
\Gamma_1(z) & =\frac{\sqrt{1-\alpha^2}}{1-\alpha z^{-1}} \\
\Gamma_2(z) & =\frac{\sqrt{1-\alpha^2}}{1-\alpha z^{-1}}\frac{z^{-1}-\alpha}{1-\alpha z^{-1}} \\
& \vdots \\
\Gamma_N(z) & =\frac{\sqrt{1-\alpha^2}}{1-\alpha z^{-1}}(\frac{z^{-1}-\alpha}{1-\alpha z^{-1}})^{N-1} \\
\end{split}
\end{equation}
where $0 \leqslant \alpha \leqslant 1$ is the pole of the discretized Laguerre network called scaling factor and $N$ is the number of terms. The orthonormality of the Laguerre functions in the time domain is expressed by:
\begin{align}
\label{19}
\sum_{k=0}^{\infty}l_i(k)l_j(k)&=0 \qquad for \quad i\neq j.\\
\label{20}
\sum_{k=0}^{\infty}l_i(k)l_j(k)&=1 \qquad for \quad i = j.
\end{align}
The discrete Laguerre function expression is simplified into a vector form as the following:
\begin{equation}
    \label{eq21}
L(k) = [l1(k)\ l_2(k)\hdots l_N(k)]^T
\end{equation}
where $l_i(k)$ is the inverse z-transform of $\Gamma_i(z,\alpha)$:  
\begin{equation}
\label{eq22}
\Gamma_k(z)=\Gamma_{k-1}(z)\frac{z^{-1}-\alpha}{1-\alpha z^{-1}}
\end{equation} 
Based on equation (\ref{eq22}), the Laguerre sequence can be expressed recursively by:
\begin{equation}
\label{eq23}
L(k+1) = A_l L(k)
\end{equation} 
with $A_l$ being an $N\times N$ matrix and a function of $\alpha$ and $\beta = (1-\alpha^2)$, and $L(0)$ as initial condition given by: 
$$L(0)^T=\sqrt{\beta}[1\ -\alpha\quad  \alpha^2\ -\alpha^3\ \hdots\  (-1)^{N-1} \alpha^{N-1}].$$
For N=5, we have:
$A_l = \left[\begin{array}{ccccc} \alpha&0&0&0&0\\
\beta&\alpha&0&0&0\\
-\alpha\beta&\beta&\alpha&0&0\\
\alpha^2\beta&-\alpha\beta&\beta&\alpha&0\\
-\alpha^3\beta&\alpha^2\beta&-\alpha\beta&\beta&\alpha \end{array}\right] $,
$L(0)^T=\sqrt{\beta} \left[\begin{array}{c} 1\\-\alpha\\ \alpha^2 \\
-\alpha^3 \\
 \alpha^4 \end{array}\right].$\\
\noindent
For the case of $\alpha = 0$ the Laguerre function becomes just a set of pulses and the MPC design becomes equivalent to the traditional approach discussed earlier. Thus, the control increment at $(k+i)$ is represented by the Laguerre terms as the following: 
\begin{equation}\label{eq24}
\Delta u(k+i)= \sum_{j=1}^{N}c_j(k)l_j(i) = L(i)^T \eta,\ i=1,2,...,N_c
\end{equation}
where $\eta =[c_1\ c_2\ \hdots\ c_N]^T$ is the Laguerre fitting coefficient, solving the QP problem becomes a matter of finding the optimal fitting coefficient $\eta_{opt}$. The solution of the QP problem given in equations (\ref{eq17}) can be divided into two cases; minimization with equality constraints, and minimization with inequality constraints. For equality constraints the QP problem is solved using Lagrange multipliers:
\begin{equation}\label{eq25}
J=\frac{1}{2}x^TEx+x^T+\lambda^T(Mx-\gamma)
\end{equation}
The Lagrange expression is subject to equality constraints ($Mx = \gamma$) and the solution is obtained by finding the optimal Lagrange multiplier $\lambda$ and corresponding decision variable $x$ through partial derivatives equated to zero:
\begin{align}
\frac{\partial J}{\partial x} & = Ex+F+M^T \lambda = 0 \label{26} \\ 
\frac{\partial J}{\partial \lambda} & = Mx-\gamma = 0 \label{27}
\end{align}
Hence, the solution is the optimal $\lambda_{opt}$ which contains the Lagrange multipliers, and the optimal decision variable $x_{opt}$ which denotes $\Delta U$ in the context of predictive control:
\begin{align}
\lambda_{opt} & = -(ME^{-1}M^T)^{-1}(\gamma+ME^{-1}F) \label{28}\\
x_{opt} & = -E^{-1}(M^T\lambda_{opt}+F) \label{29}
\end{align}
For inequality constraints, $(Mx \leqslant \gamma)$ is divided into active $(Mx=\gamma)$ and inactive $(Mx < \gamma)$ constraints and these are defined in terms of Lagrange multipliers through the Kuhn-Tucker conditions:
\begin{equation}\label{30}
\begin{split}
Ex+F+M^T\lambda& = 0\\
Mx -\gamma & \leqslant 0\\
\lambda^T(Mx-\gamma)& = 0\\
\lambda & \geqslant 0
\end{split}
\end{equation}  
The solution in this case requires the identification of active constraints first and then solving the optimization for those active constraints while inactive constraints are eliminated. In this paper, the Hildreth's QP method has been used to solve the QP problem of (\ref{eq25}), for ample details on this method we refer the reader to \cite{18}.

\section{Controller Optimization With Improved PSO}

The particle swarm optimization is a swarm intelligence algorithm that mimics animal group behaviours, such as bird flocks and fish swarms. It is an evolutionary population-based algorithm where each individual (particle) represents a possible solution and the group of individuals make up a swarm. Due to its flexibility and ease of implementation, the PSO tool is used in various fields \cite{19,20,21}. In an optimization problem with N dimensions, each $i^{th}$ particle is attributed with a velocity vector $v_i = [v_{i1},v_{i2},...,v_{in}]$ and a position vector $p_i=[p_{i1},p_{i2},...,p_{in}]$. The classic algorithm is defined by:
\begin{equation}
\label{eq31}
\left\{
    \begin{array}{ll}
 v_i(k+1) =& \omega v_i(k) + c_1 r_1 (Pb_i(k)-x_i(k)) \\
 &+ c_2 r_2(Gb(k)-x_i(k))\\
 x_i(k+1)  =& x_i(k) + v_i(k+1)
    \end{array}
\right.
\end{equation}
where at each iteration $k$, $v_i$ and $p_i$ are the velocity and position of particle $i$, and $\omega$ is the inertia weight. The coefficients $c_{1,2}$ are learning factors called cognitive and social accelerations respectively, and $r_{1,2} \in [0,1] $ are random numbers. $Pb_i$ represents the position with the best fitness score for particle $i$ and $Gb$ is the global best position in the swarm. In the classic algorithm $\omega$ and $c_{1,2}$ are constant values, several improved versions have been reported in the literature such as \cite{22} and \cite{23} where $w$ and $c_{1,2}$ are dynamic according to equation (\ref{eq32}) and table (\ref{tab:1}):

\begin{equation}
\label{eq32}
\omega = \omega_{max} - (\omega_{max}-\omega_{min}) \frac{g}{G}
\end{equation}
where $\omega_{max}$ and $\omega_{min}$ are the maximum and minimum of the inertia weight, $g$ and $G$ represent the number of the actual generation and the maximum number of generations respectively. 

\begin{table}[htb]
\caption{Rules for updating $c_1$ and $c_2$.} 
\label{tab:1}
\centering
\begin{tabular}{c c c} 
\hline
Phase & $c_1$ & $c_2$\\[0.8ex] 
\hline
Exploration & Big increase & Big decrease \\[0.8ex] 
Exploitation &  Small increase & Small decrease \\[0.8ex] 
Convergence & Small increase & Small increase \\ [0.8ex] 
Jumping out & Big decrease & Big Increase \\ [0.8ex] 
\hline 
\end{tabular}
\end{table}

In table (\ref{tab:1}), the cognitive acceleration coefficient $c_1$ is increased during the exploration phase to explore the swarm and pull the particle towards its best position $Pb$. During the exploitation phase, the social acceleration coefficient $c_2$ is increased to push the swarm towards the global best $Gb$. In the convergence phase, the swarm begins to find the globally optimal region and increasing $c_2$ while slightly decreasing $c_1$ helps accelerate the convergence. At the finale phase of the optimization, the global best particle moves from a local optimum, where most particles are clustered, towards a better optimum (the new global best). Therefore, all clustered particles must move as fast as possible towards the new global best. To achieve this, $c_2$ is greatly increased during this phase while $c_1$ is greatly decreased. 
to balance the four phases in table (\ref{tab:1}), we control $c_1$ and $c_2$ by the following logic which promotes the global search capabilities of PSO: 
\begin{equation}
\label{eq34}
\left\{
    \begin{array}{ll}
 &c_1(k+1) = c_1(k) + \alpha \\
 &c_2(k+1) = c_2(k) + \beta \\
 &\alpha = -\beta = 0.05\quad \text{for} \quad \frac gG \leq 20\%\\
 &\alpha = -\beta = 0.02\quad \text{for} \quad 20\% \leq \frac gG \leq 35\%\\
  &\alpha = -\beta = -0.035\quad \text{for} \quad 35\% \leq \frac gG \leq 75\%\\
 &\alpha = -\beta = -0.0015\quad \text{for} \quad \frac gG \geq 75\%
    \end{array}
\right.
\end{equation}
where $g$ is the current generation and $G$ is the maximum number of generations.\\

According to (\ref{eq32}), $\omega$ decreases linearly over iterations. A small $\omega$ promotes the local search of PSO \cite{23}. A new formula for updating the inertia weight is proposed in this paper (\ref{eq33}), which ensures an exponential decrease rather than linear. 

\begin{equation}
\label{eq33}
\omega = \omega_{min} + \frac{\exp{(\omega_{max}-\lambda_1(\omega_{max}+\omega_{min})\frac{g}{G})}}{\lambda_2}
\end{equation}
This has been found to effectively enhance the overall search capabilities of PSO when coupled with the control logic of $c_{1,2}$. Unlike a linear decrease, an exponential decrease of $\omega$ ensures faster convergence and better overall optimization. $\lambda_{1,2}$ are adjustable constants to achieve a decrease from $\omega_{max}$ to $\omega_{min}$.

\begin{figure}[t]
\centering
\includegraphics[width=0.4\textwidth]{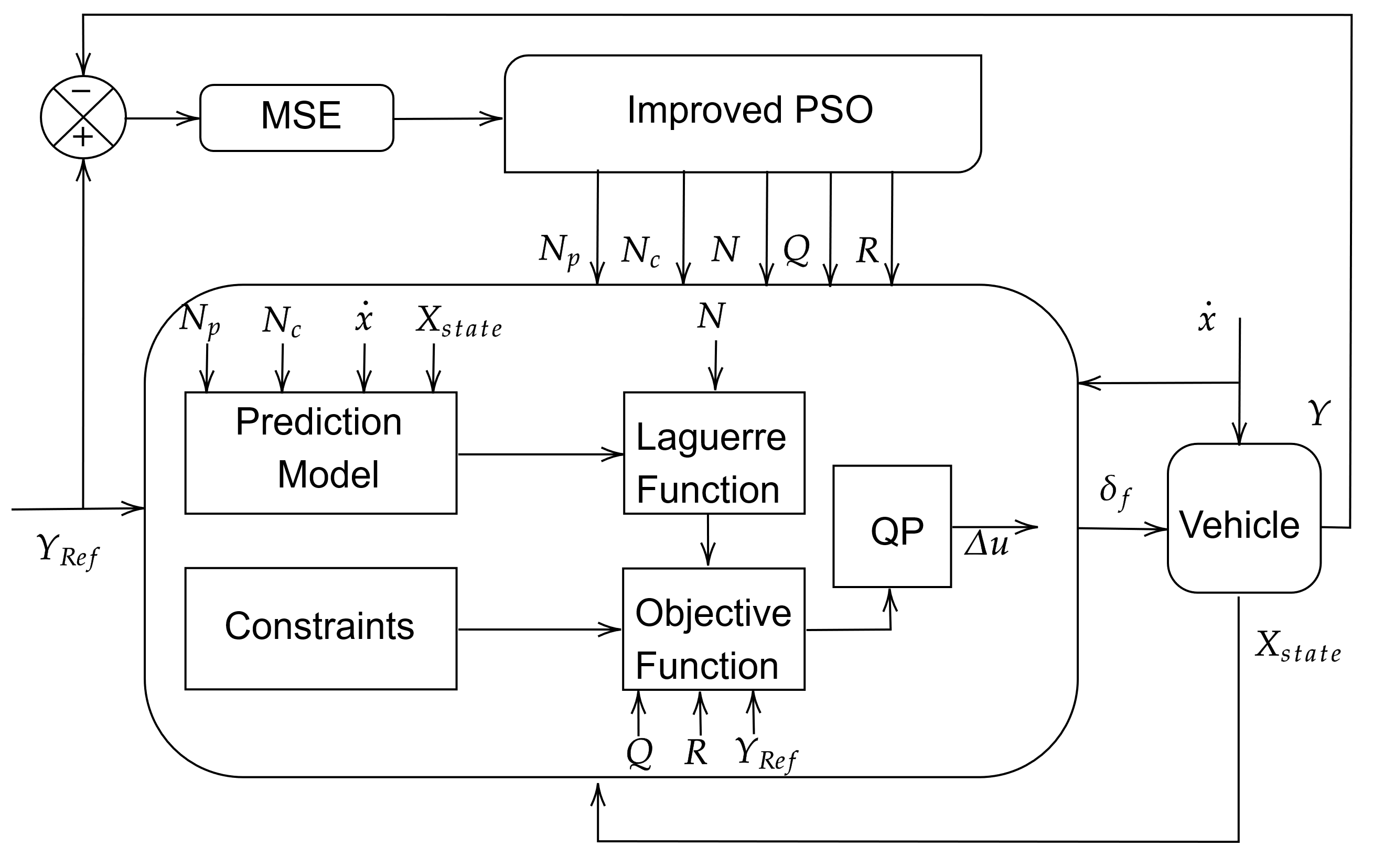}
\caption{Schematic diagram  of the proposed approach.}
\label{fig:3}
\end{figure}

The proposed PSO algorithm is used to optimize the parameters of the MPC designed with the Laguerre function in section \ref{sec:mpc}. These parameters include; the prediction and control horizons $N_p$ and $N_c$, the weighting matrices $Q$ and $R$ and the number of Laguerre terms $N$. To the best of our knowledge, existing research has not addressed the tuning and adaptation of these parameters all-together with the proposed PSO algorithm. Thus, the proposed approach allows us to optimize the controller performance in terms of its tunable parameters and Laguerre approximation at the same time without losing the computation efficiency. The mean squared error (MSE) is used as a fitness function to assess the optimality of the solutions. The schematic diagram of our approach is given in Fig. \ref{fig:3}, where $Y_{ref}$, $Y$, $\dot{x}$ and $X_{state}$ are the reference path, the vehicle lateral position, the longitudinal velocity and the states vector respectively.  
\setlength{\textfloatsep}{2pt}

\section{Results and Discussion}

\begin{figure}[t]
\centering
\includegraphics[width=8cm,height=3.5cm]{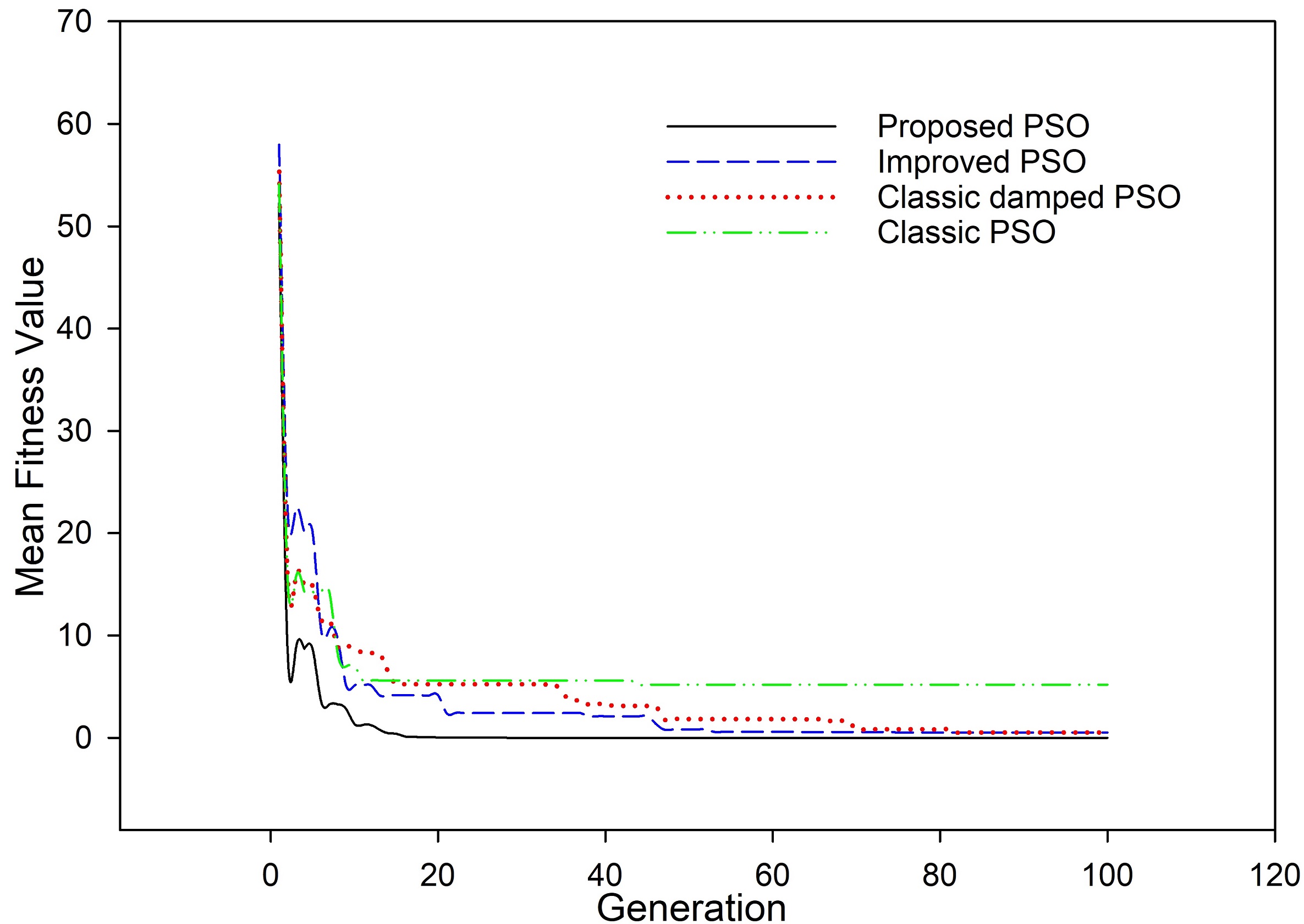}
\caption{Performance of the proposed PSO.}
\label{fig:4}
\end{figure}

The performance of the proposed PSO algorithm is compared to the improved version of \cite{23}, the modified PSO of \cite{22} and the latter enhanced with a damped inertia weight instead of a constant one ($\omega(k+1) = \omega(k)\times \omega_{dmp}$). The sphere function ($f_1(x) = \sum_{i=1}^nx_i^2$) is used as a benchmark test similar to \cite{23}, and the optimization is carried out with the parameters given in table (\ref{tab:2}) and the logic stated by equations (\ref{eq34}). The obtained results in Fig. \ref{fig:4}, show that the proposed PSO is much faster than the other PSO versions and able to find better optimal solutions. After $41$ generations, the proposed PSO reached a fitness of $10^{-3}$ compared to $2.1$ and $3.11$ for the improved PSO and the classic damped PSO respectively, while the classic PSO reached only $5.59$.

\begin{table}[b]
\caption{PSO hyper parameters.} 
\label{tab:2}
\centering
\begin{tabular}{c c c} 
\hline
Parameter & Interpretation & Value\\[0.8ex] 
\hline
$G$ & Number of generations & $15$ \\[0.8ex] 
$N_{Pop}$ &  Number of particles & $20$ \\[0.8ex] 
$\omega_{max}$ & Maximum inertia weight & $0.99$ \\ [0.8ex] 
$\omega_{min}$ & Minimum inertia weight & $0.1$ \\ [0.8ex] 
$\omega_{dmp}$ & Damping constant & $0.99$ \\ [0.8ex] 
$c_{1i}$ & Initial cognitive acceleration coefficient & $2$ \\ [0.8ex] 
$c_{2i}$ & Initial social acceleration coefficient & $2$ \\ [0.8ex] 
$\lambda_1$ & Constant & $30$ \\ [0.8ex]
$\lambda_2$ & Constant & $3$ \\ [0.8ex]
\hline 
\end{tabular}
\end{table}

\begin{figure}[!h]
\centering
\includegraphics[width=8cm,height=3.9cm]{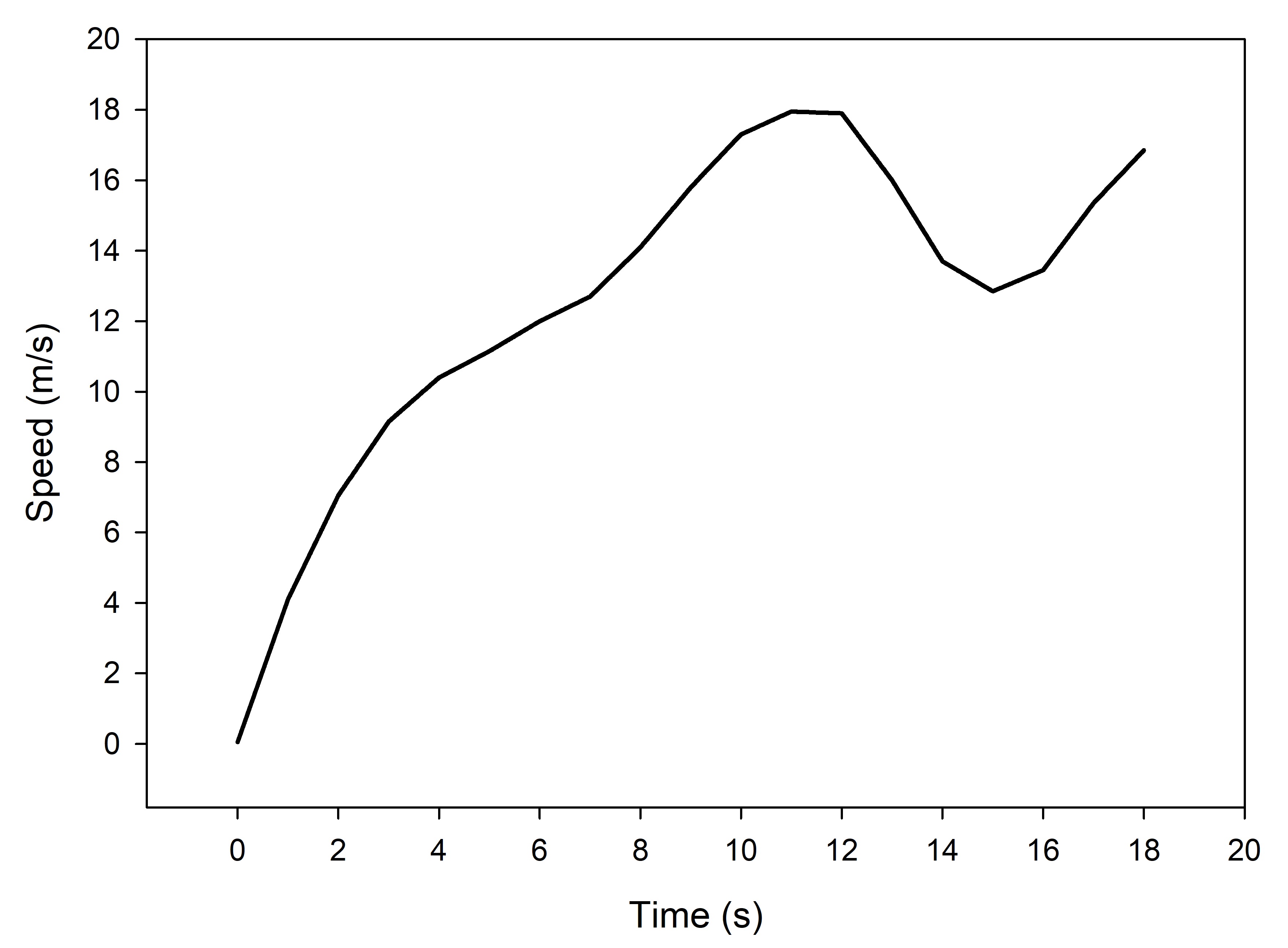}
\caption{Longitudinal velocity profile for \textbf{$Sc_2$}.}
\label{fig:8}
\end{figure}

\begin{figure}[!h]
\centering
\includegraphics[width=8cm,height=4cm]{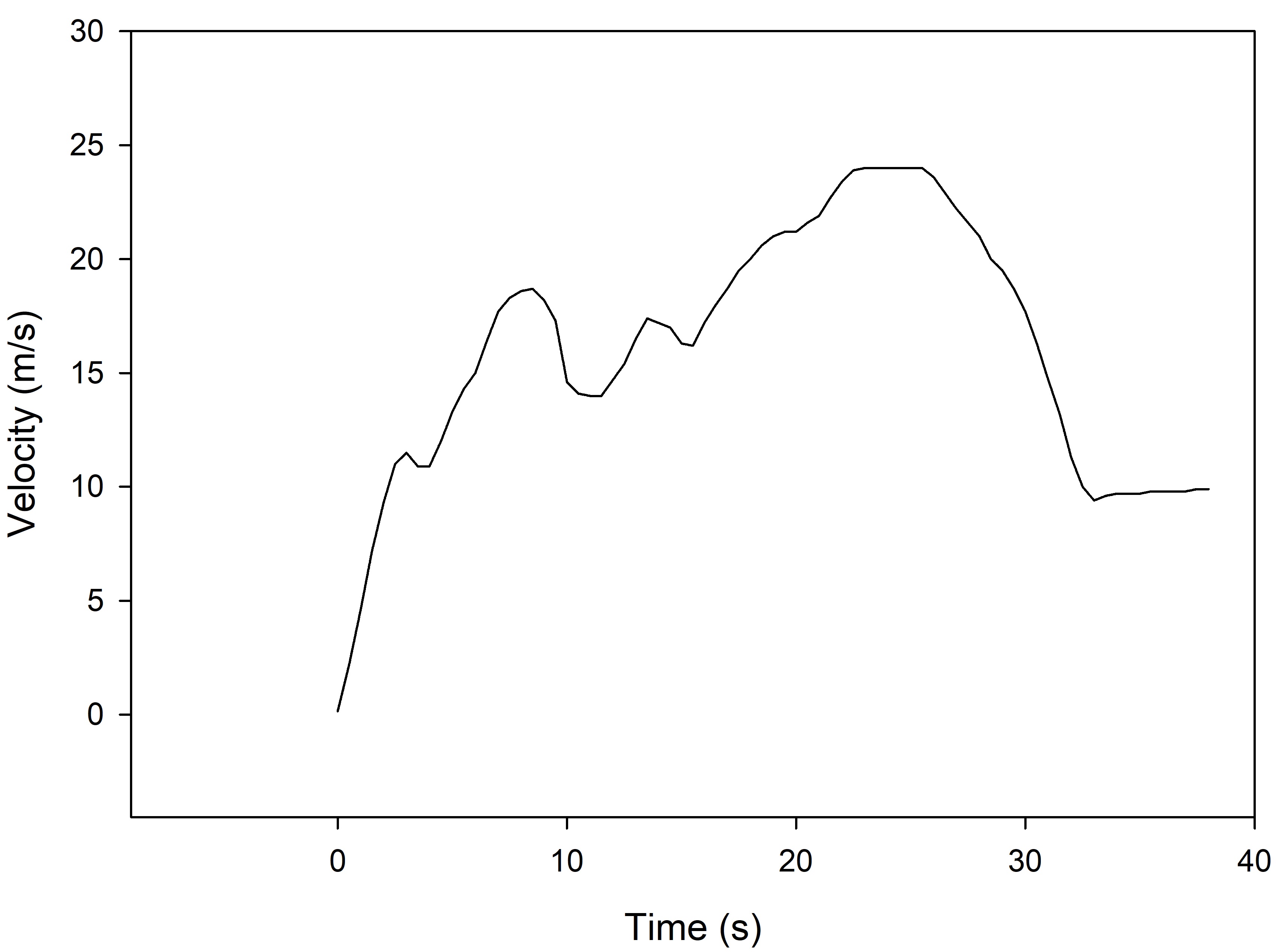}
\caption{Velocity profile for \textbf{$Sc_3$}.}
\label{fig:15}
\end{figure}

\begin{figure}[!h]
\centering
\includegraphics[width=8cm,height=3.9cm]{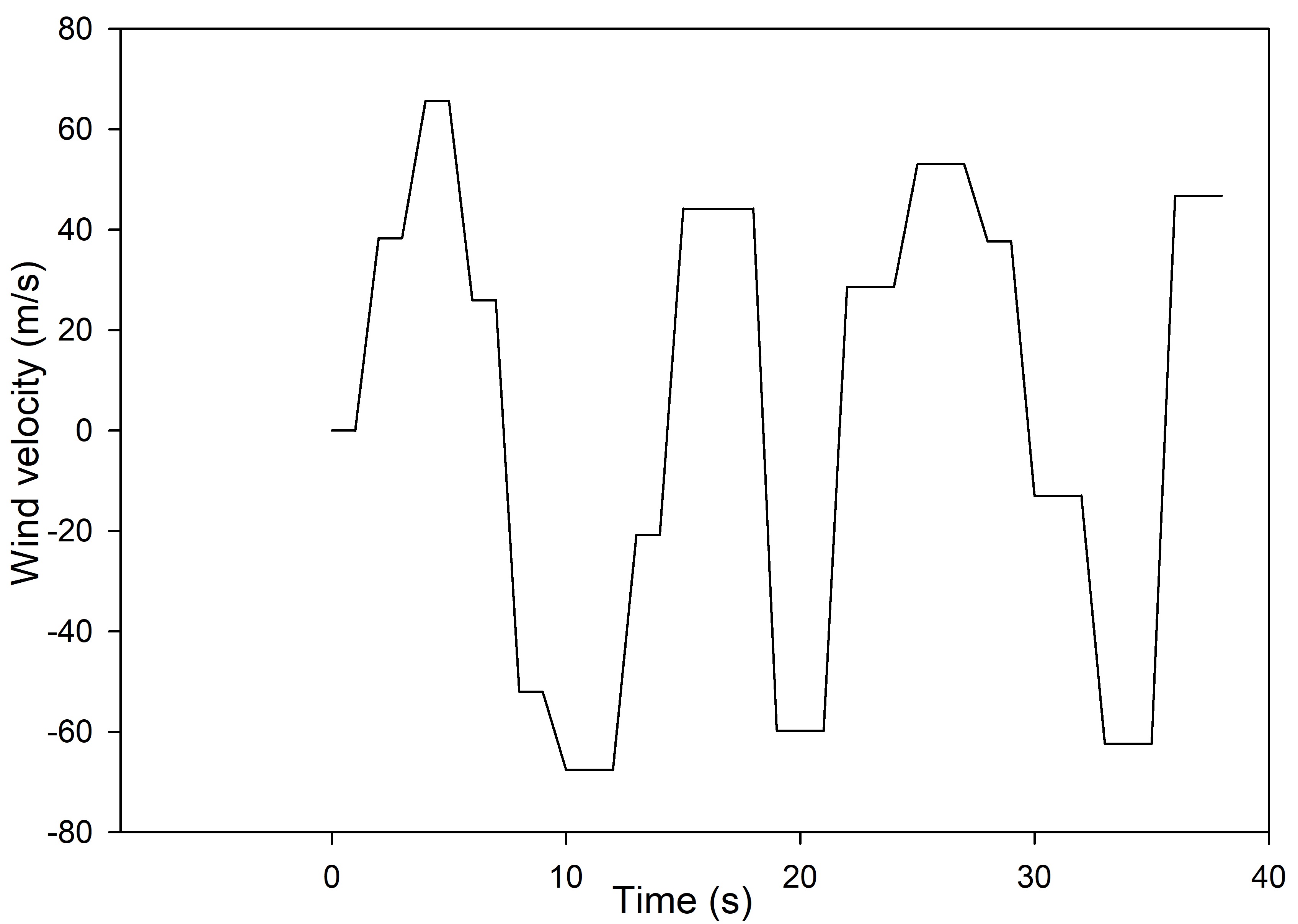}
\caption{Wind disturbance profile.}
\label{fig:21}
\end{figure}

 
A high fidelity vehicle model has been built using the Vehicle Dynamics Blockset of \texttt{MATLAB}, the latter consists of a 3 DOF dual track lateral dynamics block with nonlinear Pacejka tire formula \cite{25}. In addition, a simplified powertrain block is added to accommodate varying speed profiles with the predictive longitudinal driver block. The proposed MPC controller is optimized with the proposed PSO using $15$ generations and $20$ particles, the rest of the parameters are identical to table (\ref{tab:2}). The optimization is done for multiple longitudinal speeds ($\dot{x}(m/s) \in [3,27] $) and reference lateral positions ($y_{ref}(m)\in [-15,15]$). The adaptive MPC (AMPC) has been tested against the pure pursuit controller (P.P) and the classic MPC (tuned by trial and error) for three scenarios:
\begin{enumerate}
    \item $Sc_1$ consists of a double lane change with constant longitudinal velocity $v_x=9 m/s$. 
    \item $Sc_2$ is a double lane change with varying longitudinal velocity (Fig. \ref{fig:8}) where the AMPC is adapted to both trajectory and velocity profile.
    \item $Sc_3$ where the AMPC is tested for a general trajectory with varying longitudinal velocity (Fig. \ref{fig:15}) subject to wind gust as an external disturbance (Fig. \ref{fig:21}).
\end{enumerate}
\setlength{\textfloatsep}{2pt}

\begin{figure}[!h]
\centering
\includegraphics[width=8cm,height=3.9cm]{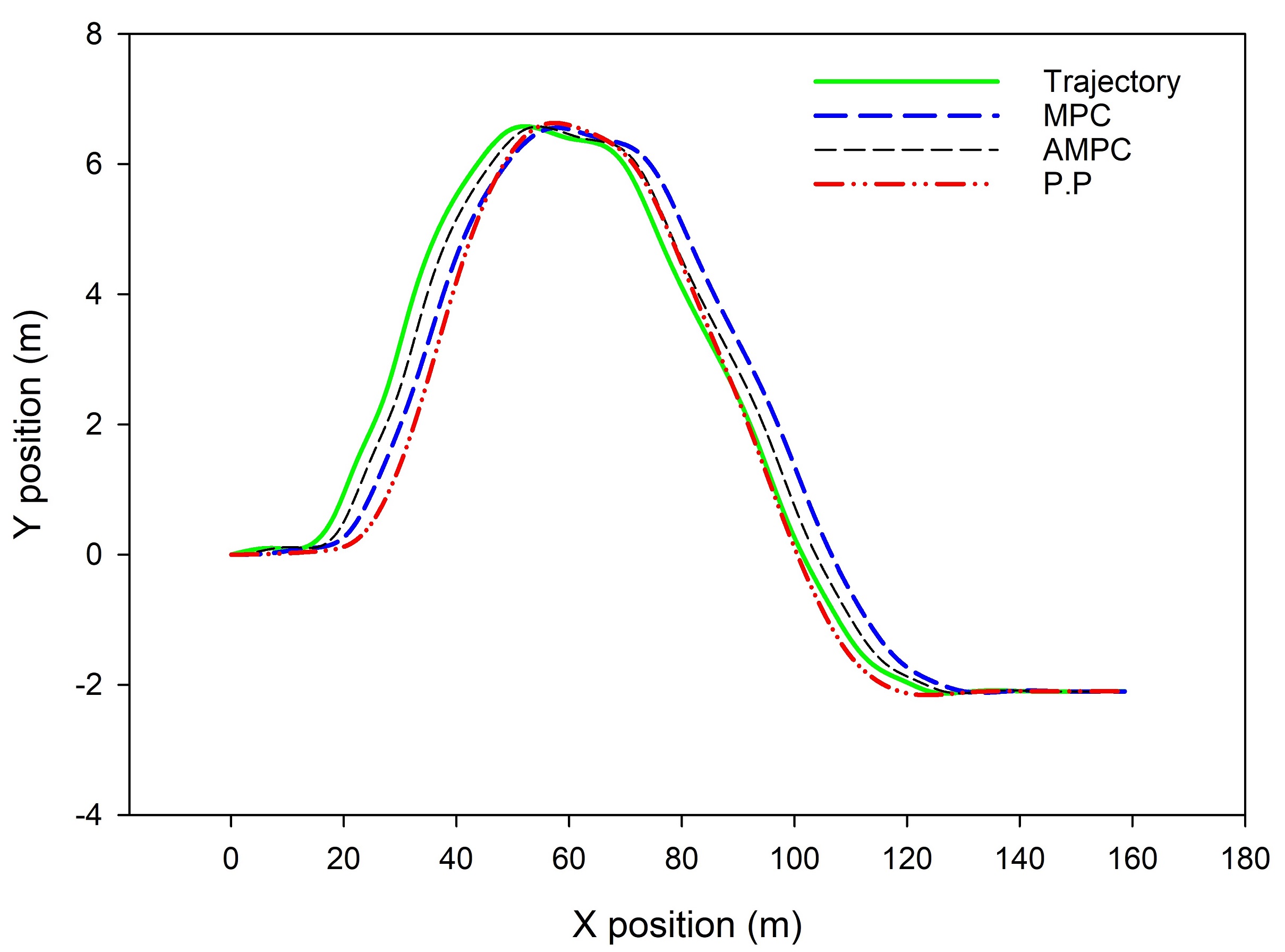}
\caption{Path tracking for \textbf{$Sc_1$}.}
\label{fig:5}
\end{figure}

\begin{figure}[!h]
\centering
\includegraphics[width=8cm,height=3.9cm]{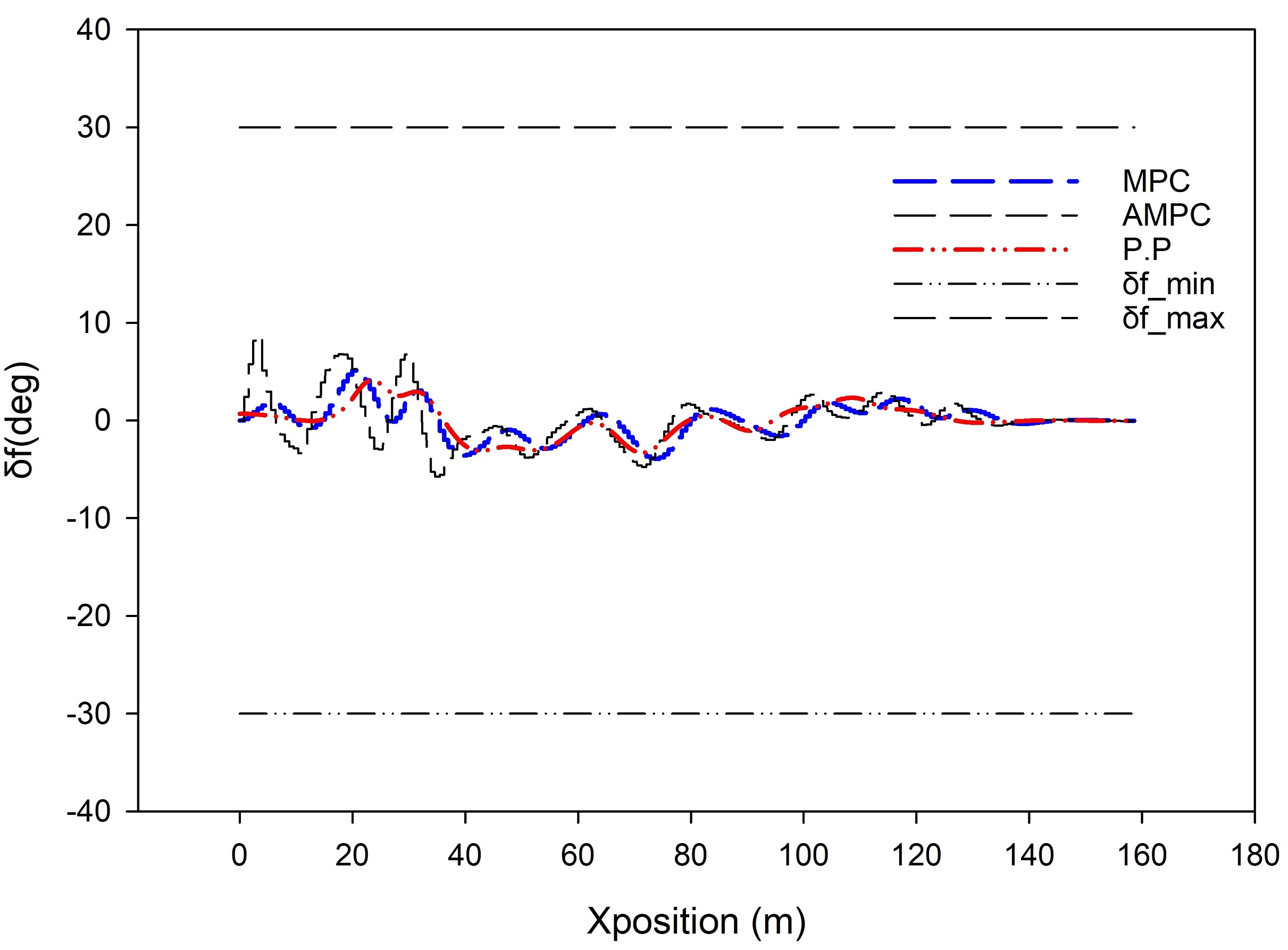}
\caption{Control signal for \textbf{$Sc_1$}.}
\label{fig:6}
\end{figure}

\begin{figure}[!h]
\centering
\includegraphics[width=8cm,height=3.9cm]{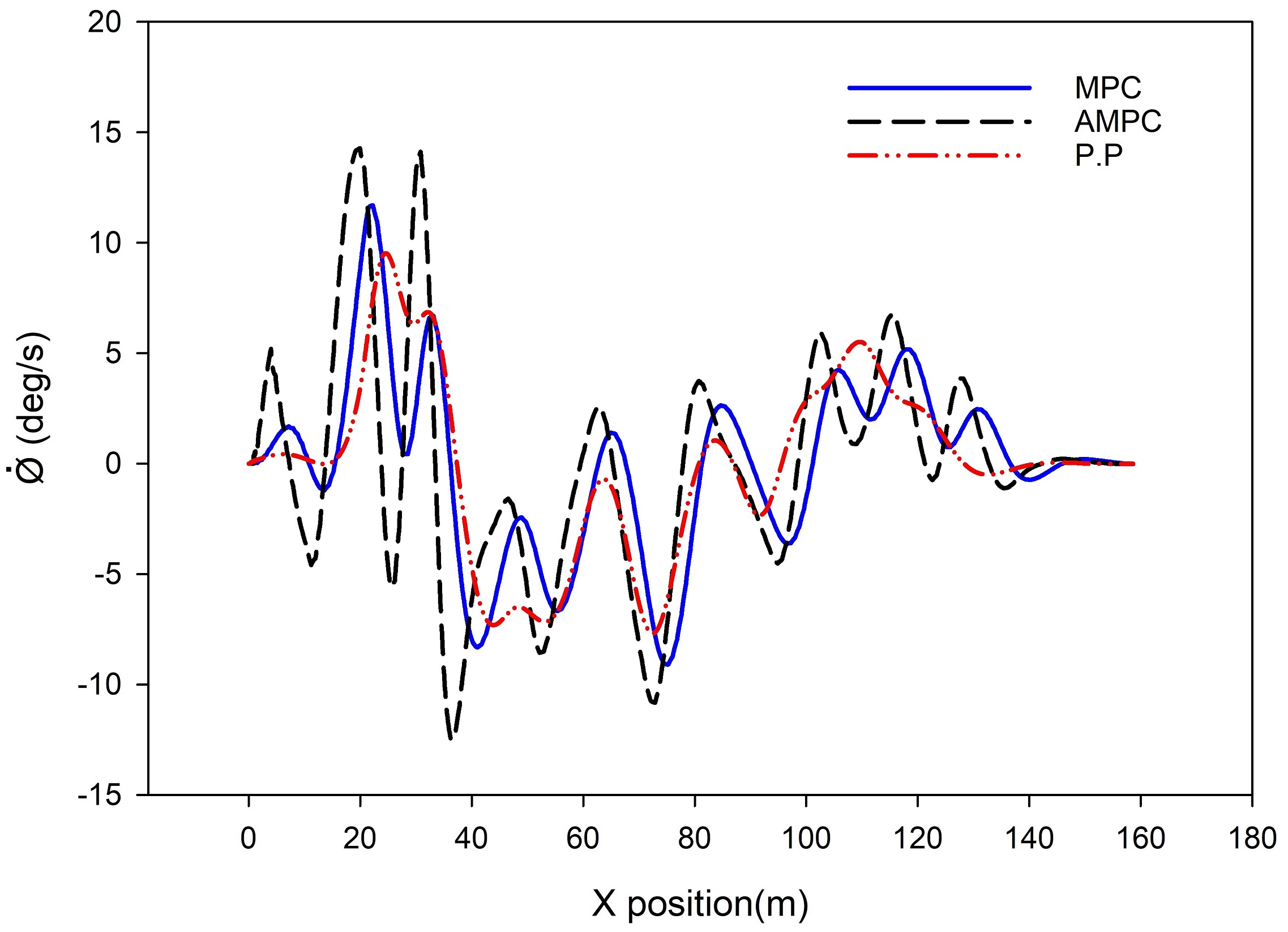}
\caption{Heading rate for \textbf{$Sc_1$}.}
\label{fig:7}
\end{figure}

The parameters of the dynamic model (\ref{eq6}) and the MPC controller are listed in table (\ref{tab:3}). The results of (\textbf{$Sc_1$}) show that the three controllers are able to track the double lane change trajectory. However, the proposed AMPC outperforms the other two controllers where it achieved an MSE value of $e_1=0.097$ compared to $e_2=0.422$ and $e_3=0.482$  for the classic MPC and pure pursuit respectively. Fig. \ref{fig:5}-\ref{fig:7} compare the tracking performance ($Y_{position}$), the steering angles ($\delta_f$) and the heading rates ($\dot{\phi}$) for the three controllers. 

\begin{table}[b]
\caption{MPC and dynamic model parameters} 
\label{tab:3}
\centering
\begin{tabular}{c c c} 
\hline
Parameter & Interpretation &Value\\[0.8ex] 
\hline
$m$ & vehicle mass & 1575 (kg) \\[0.8ex] 
$I_z$ & moment of inertia & $2875 (kg.m^2)$ \\[0.8ex] 
$l_f$ & front axle to (CG) & $1.2 (m)$ \\ [0.8ex] 
$l_r$ & rear axle to (CG) & $1.6 (m)$ \\ [0.8ex] 
$C_f$ & front cornering stiffness & $19000 (N/rad)$ \\ [0.8ex] 
$C_r$ & rear cornering stiffness & $33000 (N/rad)$ \\ [0.8ex] 
$T$ & sampling time & $0.1 (s)$ \\ [0.8ex] 
$N_p$ & Prediction horizon & $45$ \\[0.8ex] 
$N_c$ & Control horizon & $15$ \\[0.8ex] 
$\Delta u_{max/min}$ & Control increment constraints & $\pm \frac{\pi}{12}(rad)$ \\[0.8ex] 
$u_{max/min}$ & Control constraints & $\pm \frac{\pi}{6}(rad)$ \\[0.8ex] 
$R$ & Control magnitude weight & $0.01$ \\ [0.8ex] 
$Q_y$ & Tracking precision weight & $10$ \\ [0.8ex] 
$N$ & Number of Laguerre terms & $5$ \\ [0.8ex] 
$\alpha$ & Laguerre scale factor & $0.75$ \\ [0.8ex] 
\hline 
\end{tabular}
\end{table}

\setlength{\textfloatsep}{2pt}

\begin{figure}[!h]
\centering
\includegraphics[width=8cm,height=3.9cm]{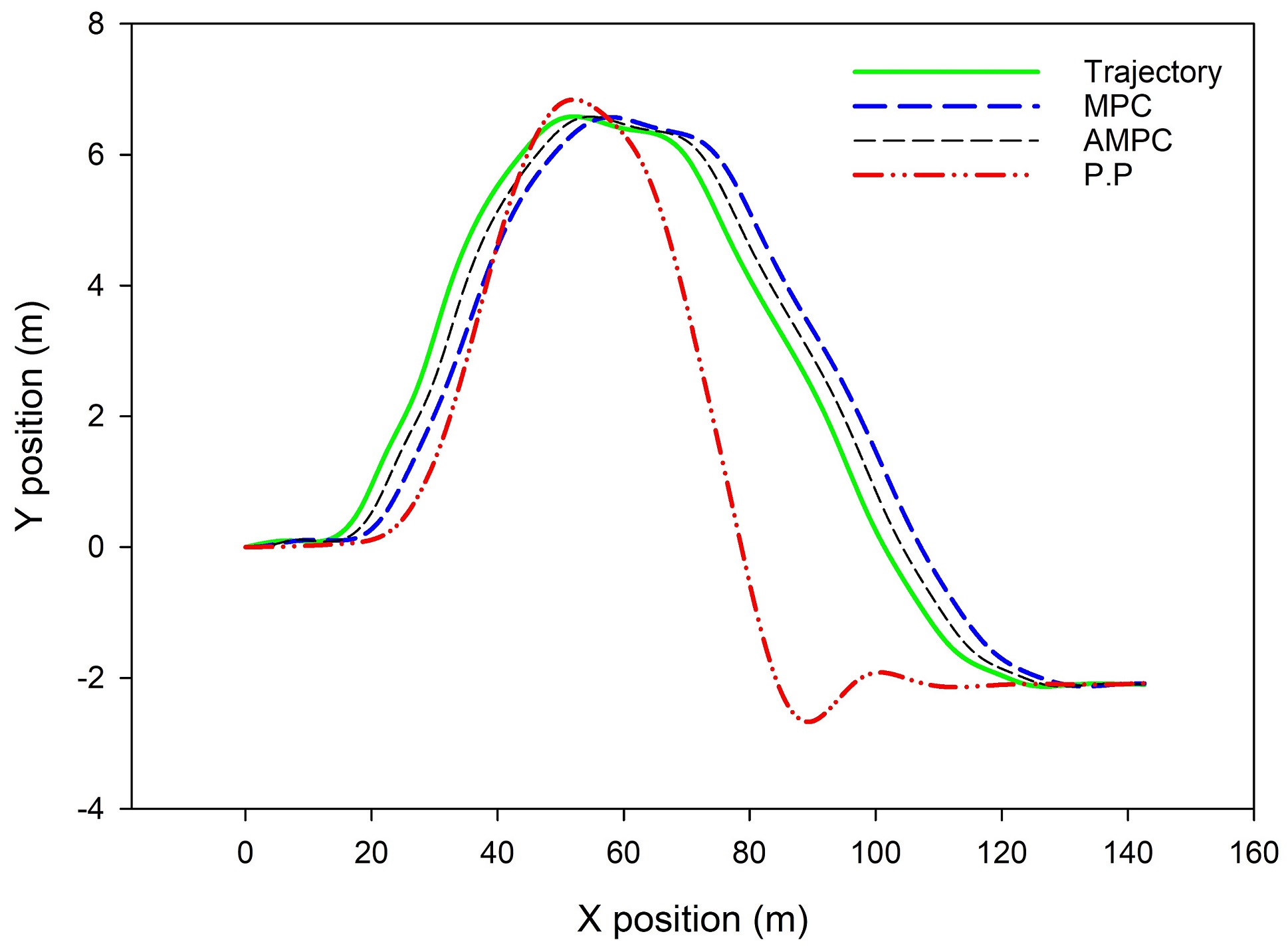}
\caption{Path tracking for \textbf{$Sc_2$}.}
\label{fig:9}
\end{figure}

\begin{figure}[!h]
\centering
\includegraphics[width=8cm,height=3.9cm]{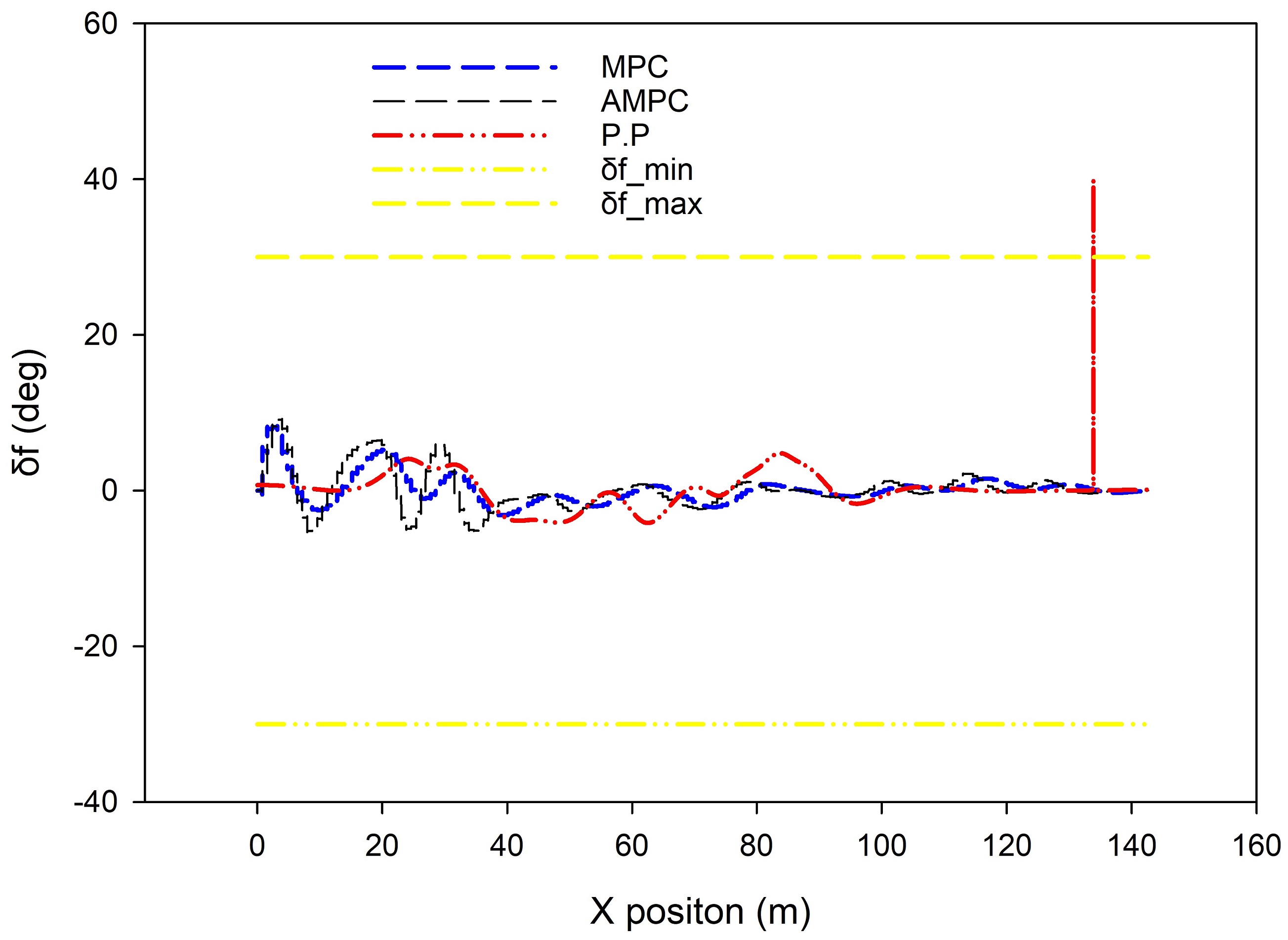}
\caption{Control signal for \textbf{$Sc_2$}.}
\label{fig:10}
\end{figure}

\begin{figure}[!h]
\centering
\includegraphics[width=8cm,height=3.9cm]{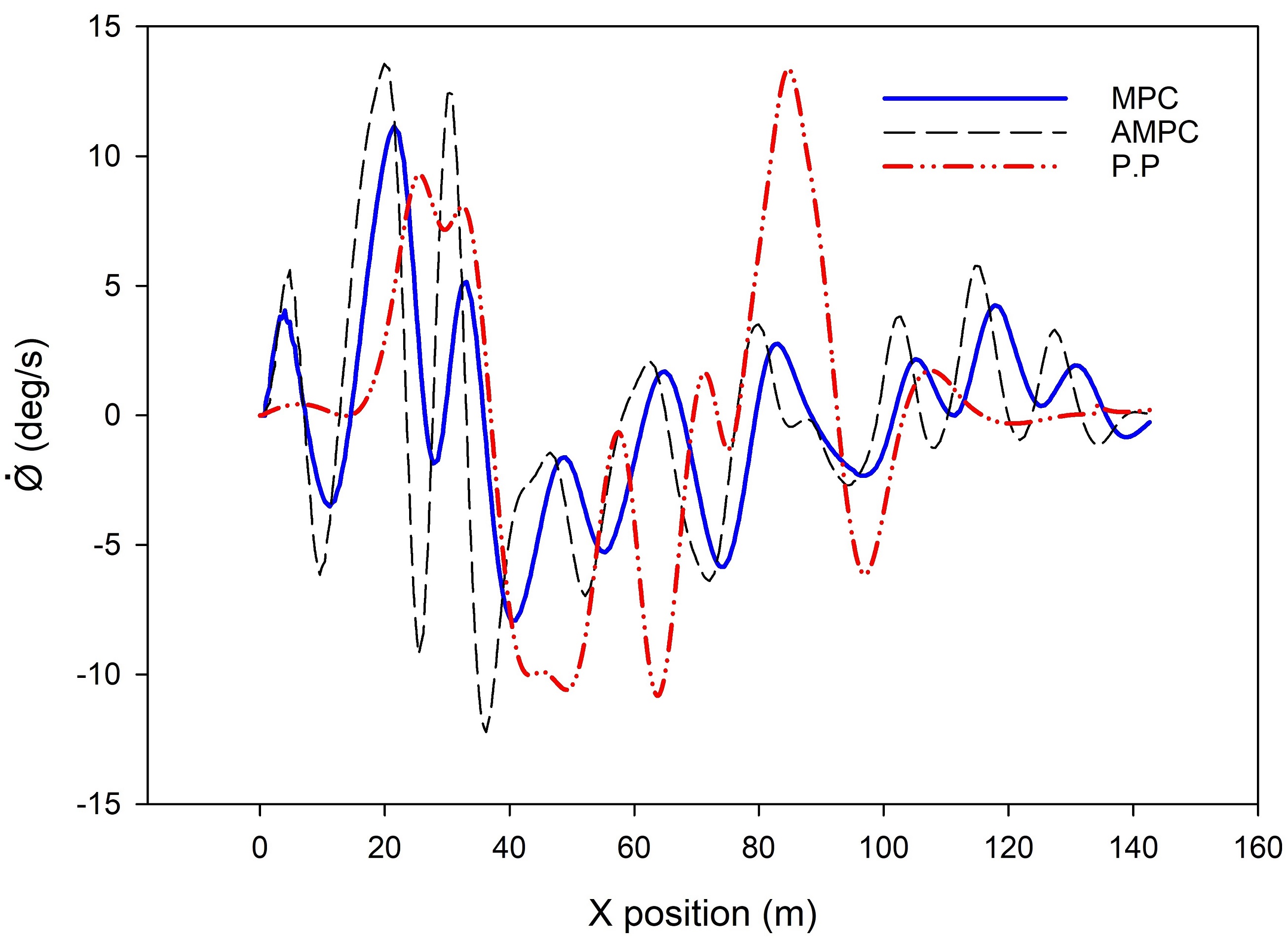}
\caption{Heading rate for \textbf{$Sc_2$}.}
\label{fig:11}
\end{figure}

On the other hand, the P.P controller could not perform the double lane change of (\textbf{$Sc_2$}), while the proposed AMPC tracked the trajectory significantly better than the classic MPC. Their respective MSE values are $e_1=0.124$ against $e_2=0.489$. Fig. \ref{fig:9}-\ref{fig:11} show the tracking performance, the control signals and the heading rates.
\setlength{\textfloatsep}{2pt}

\begin{figure}[!h]
\centering
\includegraphics[width=8cm,height=3.75cm]{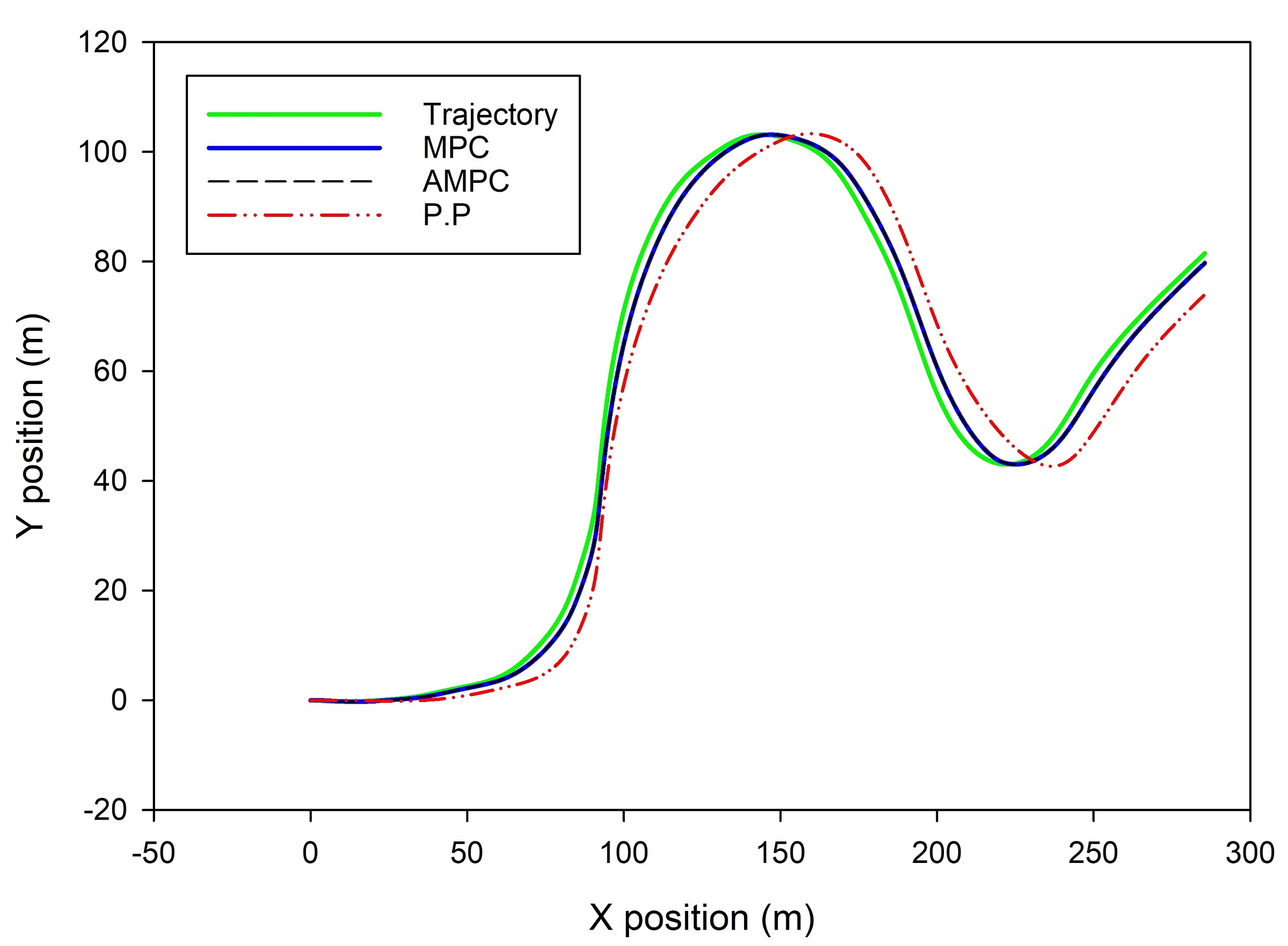}
\caption{Path tracking for \textbf{$Sc_3$} with constant velocity.}
\label{fig:12}
\end{figure}

\begin{figure}[!h]
\centering
\includegraphics[width=8cm,height=3.75cm]{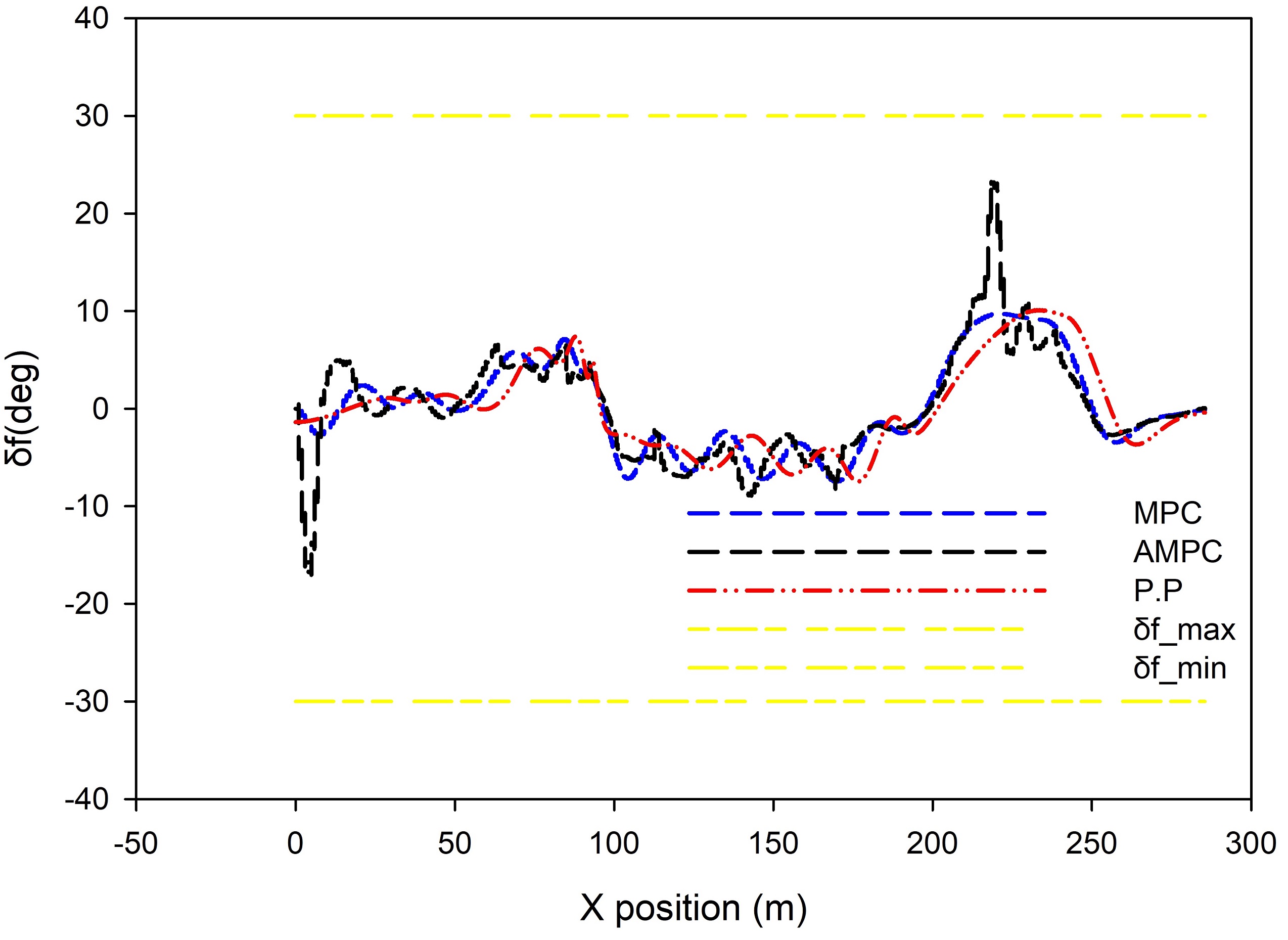}
\caption{Control signal for \textbf{$Sc_3$} with constant velocity.}
\label{fig:13}
\end{figure}

\begin{figure}[!h]
\centering
\includegraphics[width=8cm,height=3.75cm]{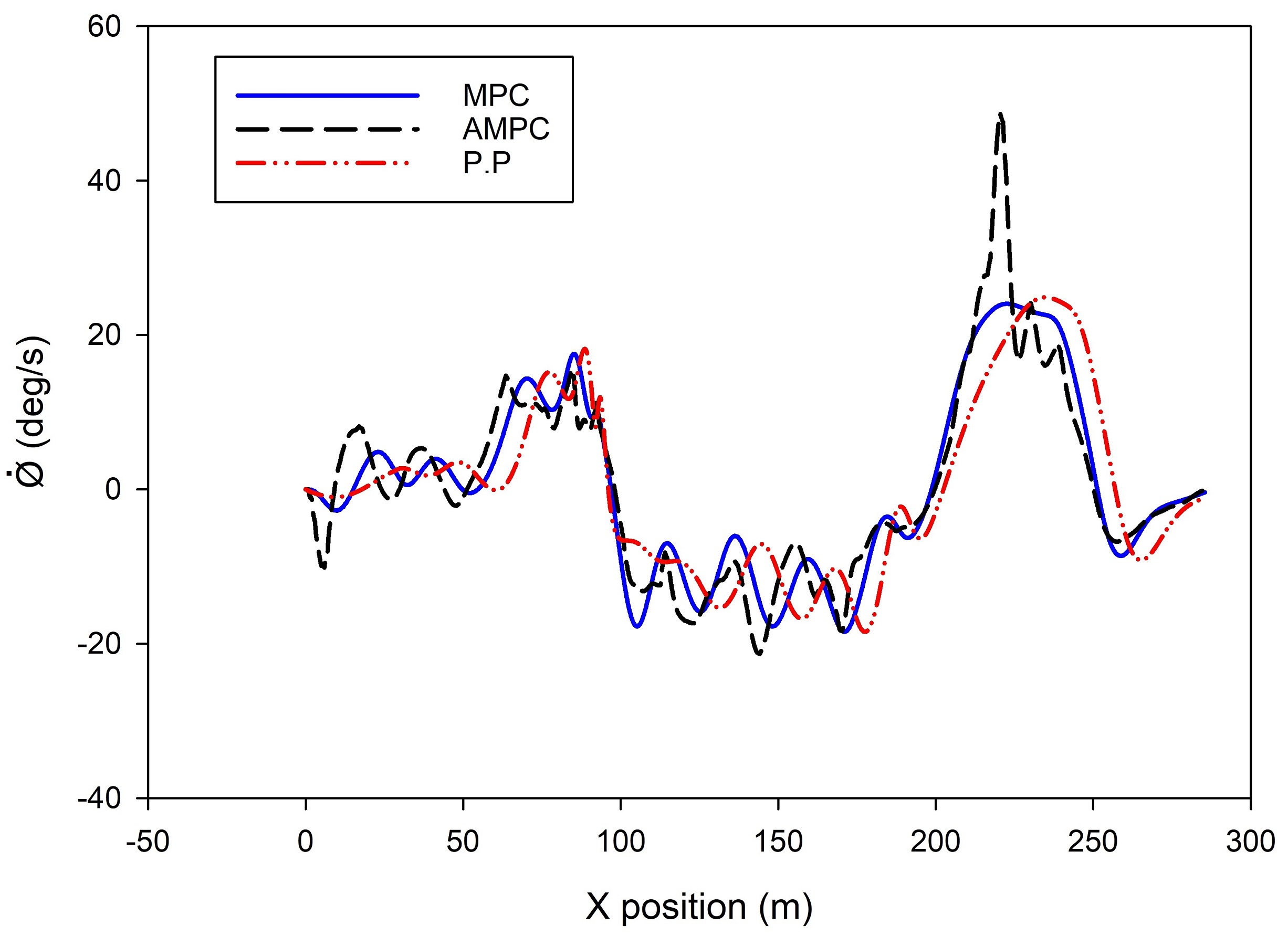}
\caption{Heading rate for \textbf{$Sc_3$} with constant velocity.}
\label{fig:14}
\end{figure}

\begin{figure}[!h]
\centering
\includegraphics[width=8cm,height=3.9cm]{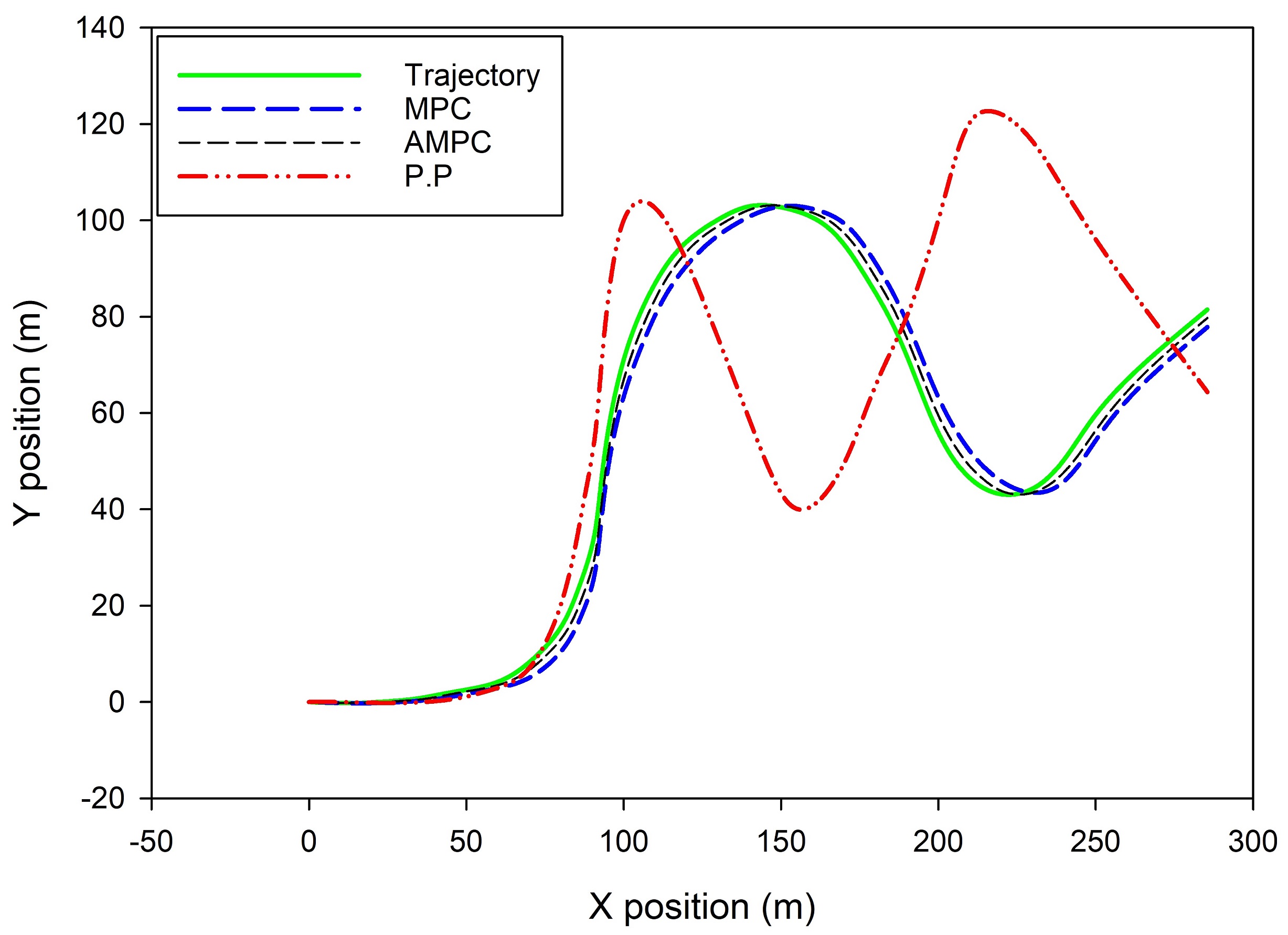}
\caption{Path tracking for \textbf{$Sc_3$} with varying velocity.}
\label{fig:16}
\end{figure}

\begin{figure}[!h]
\centering
\includegraphics[width=8cm,height=3.9cm]{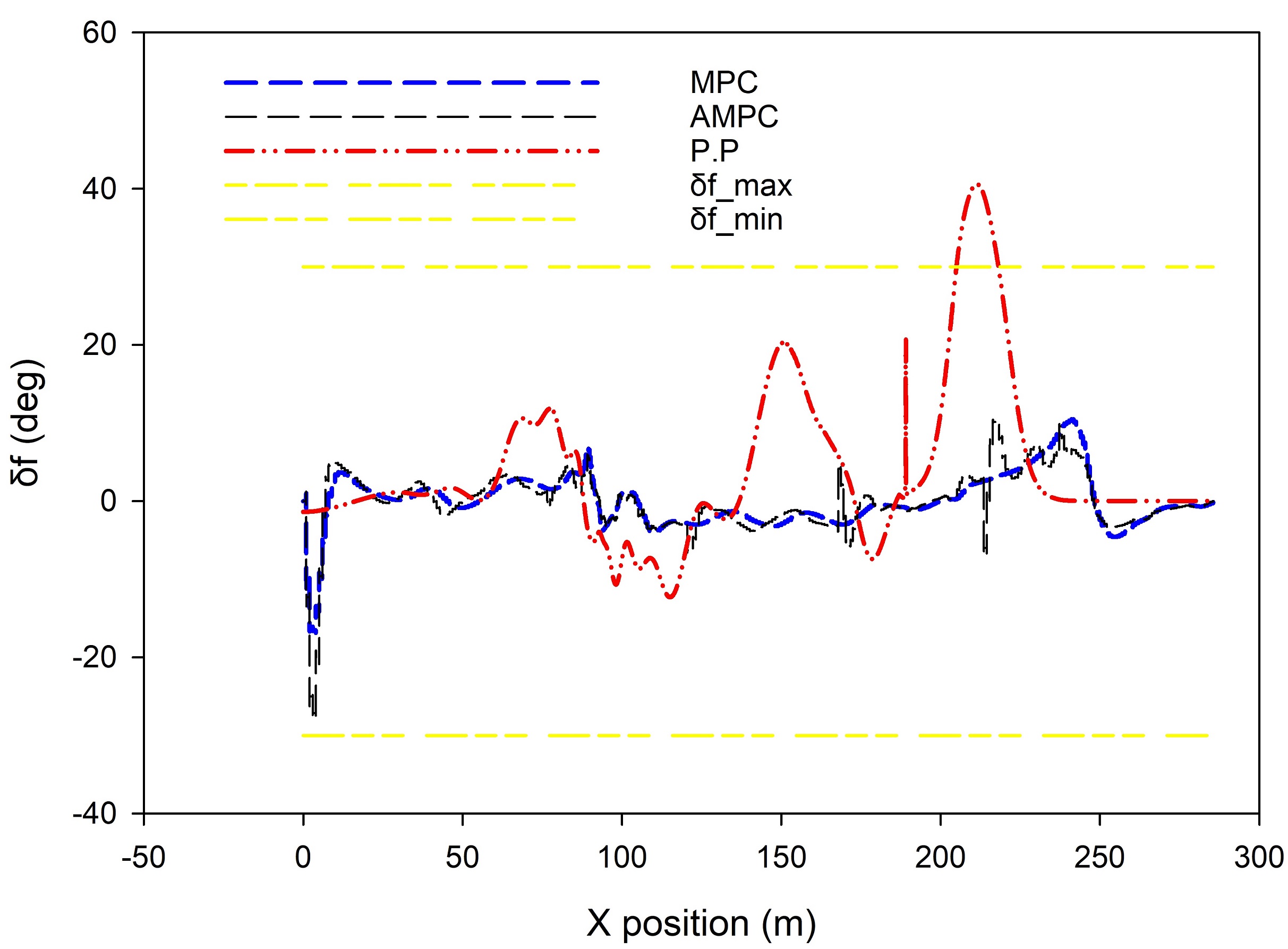}
\caption{Control signal for \textbf{$Sc_3$} with varying velocity.}
\label{fig:17}
\end{figure}

\begin{figure}[!h]
\centering
\includegraphics[width=8cm,height=3.9cm]{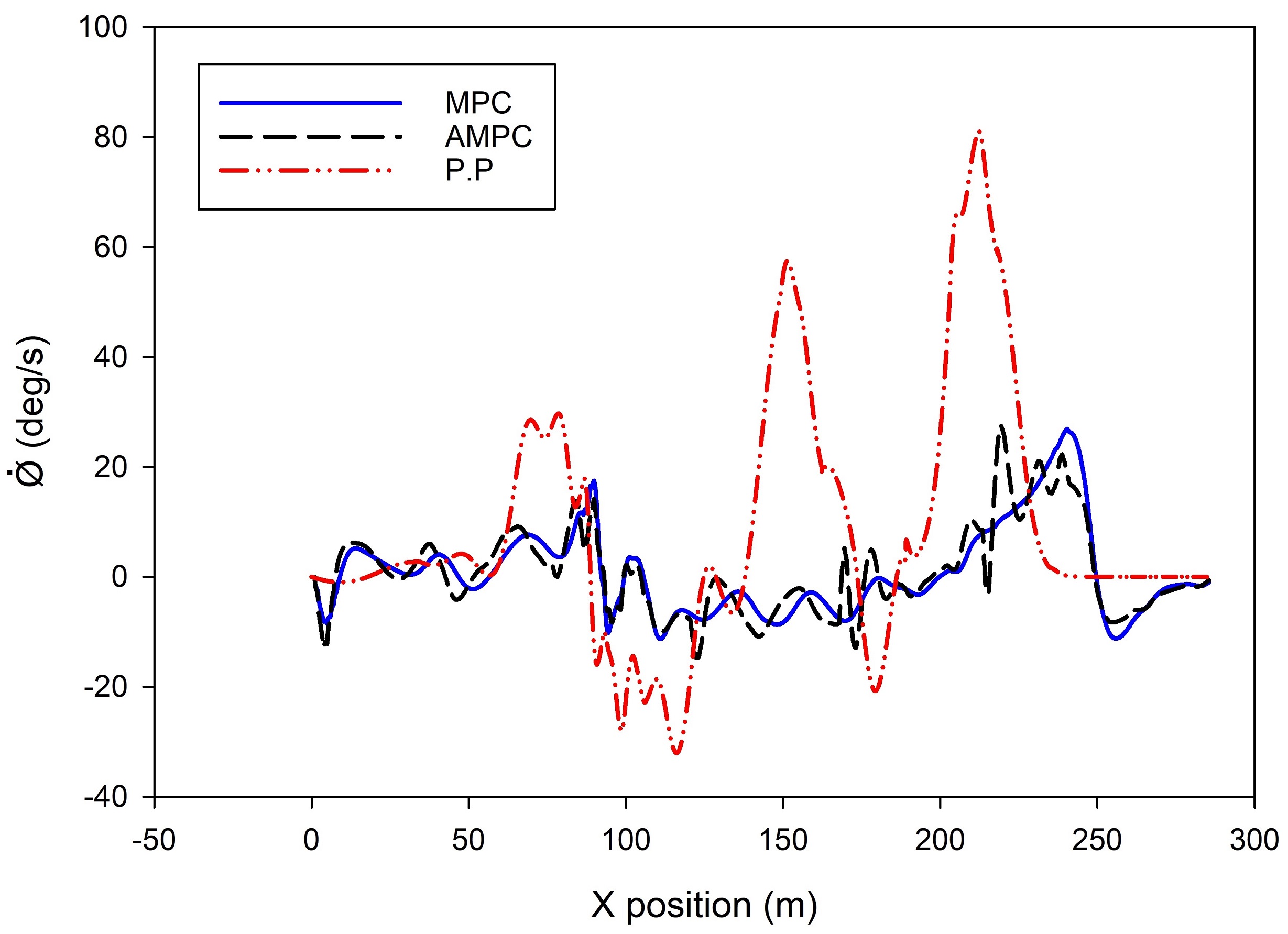}
\caption{Heading rate for \textbf{$Sc_3$} with varying velocity.}
\label{fig:18}
\end{figure}

\begin{figure}[!h]
\centering
\includegraphics[width=8cm,height=3.9cm]{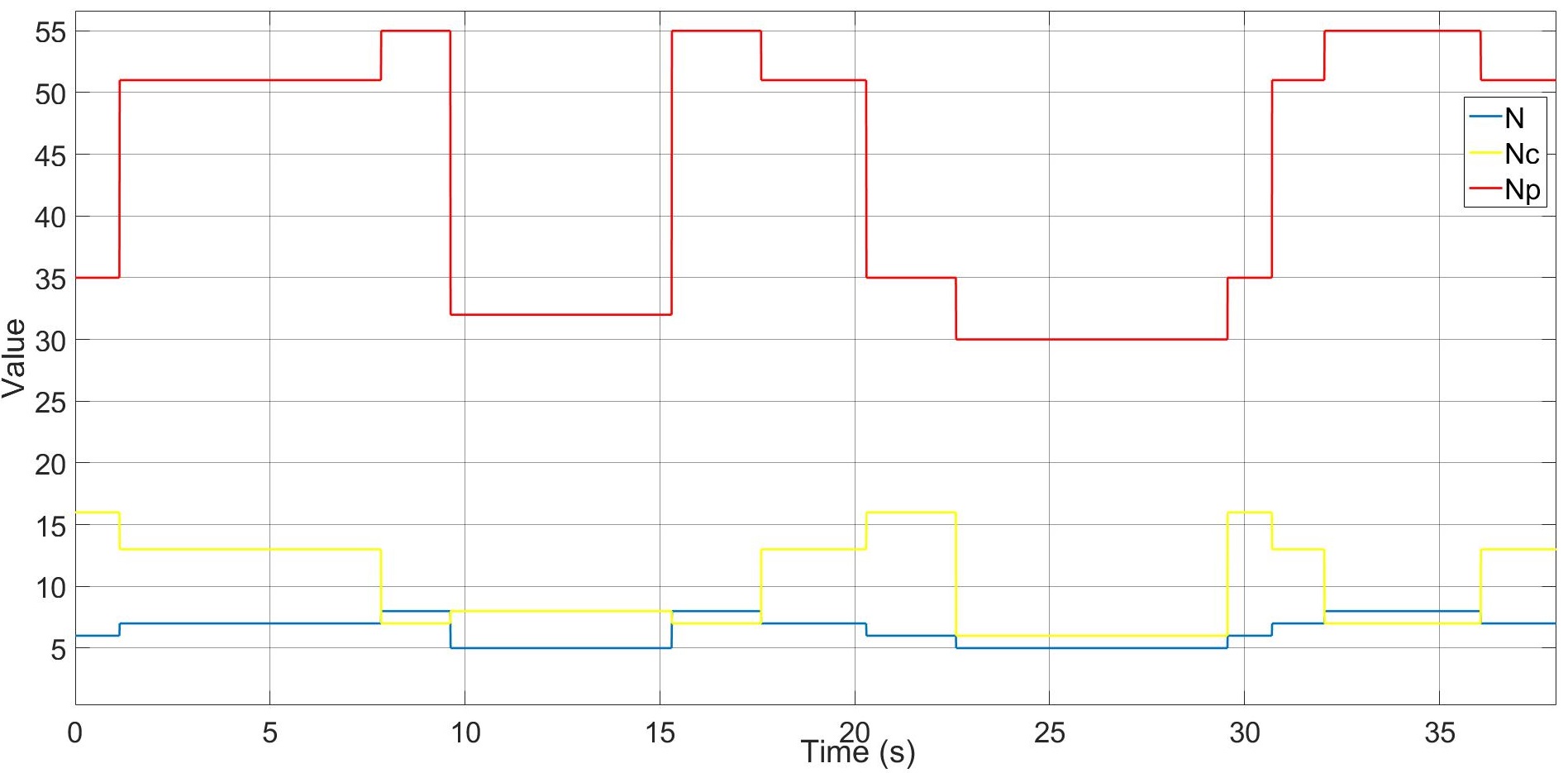}
\caption{MPC parameters adaptation for \textbf{$Sc_3$} with varying velocity.}
\label{fig:19}
\end{figure}

\begin{figure}[!h]
\centering
\includegraphics[width=8cm,height=3.9cm]{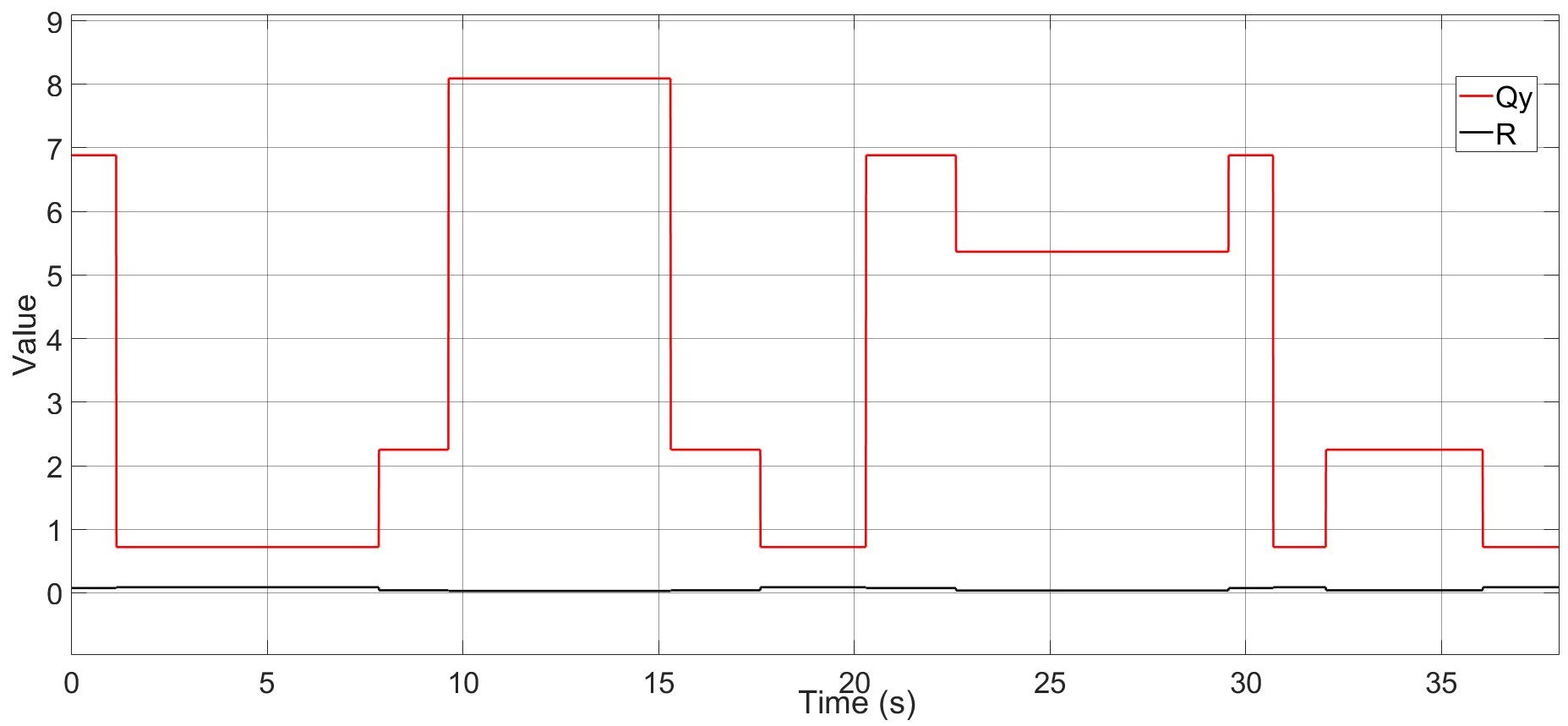}
\caption{Cost function adaptation for \textbf{$Sc_3$} with varying velocity.}
\label{fig:20}
\end{figure}

The results for (\textbf{$Sc_3$}) confirm better tracking ability of the proposed AMPC. Fig. \ref{fig:12}-\ref{fig:14} show the tracking performance for a constant velocity $v_x = 11 m/s$. The MSE values are $e_1=10.79$, $e_2=21.234$ and $e_3=76.639$ for the AMPC, MPC and P.P respectively. In Fig. \ref{fig:16}-\ref{fig:18}, the same trajectory is simulated with a varying velocity profile (Fig. \ref{fig:15}). The performance of AMPC is much improved compared to regular MPC, but both are able to track the trajectory and the resulting MSE values for AMPC and MPC are $e_1=7.219$ and $e_2=25.079$. The parameters adaptation for this case is illustrated in Fig. \ref{fig:19}-\ref{fig:20}.  
\setlength{\textfloatsep}{2pt}

\begin{figure}[!h]
\centering
\includegraphics[width=8cm,height=3.9cm]{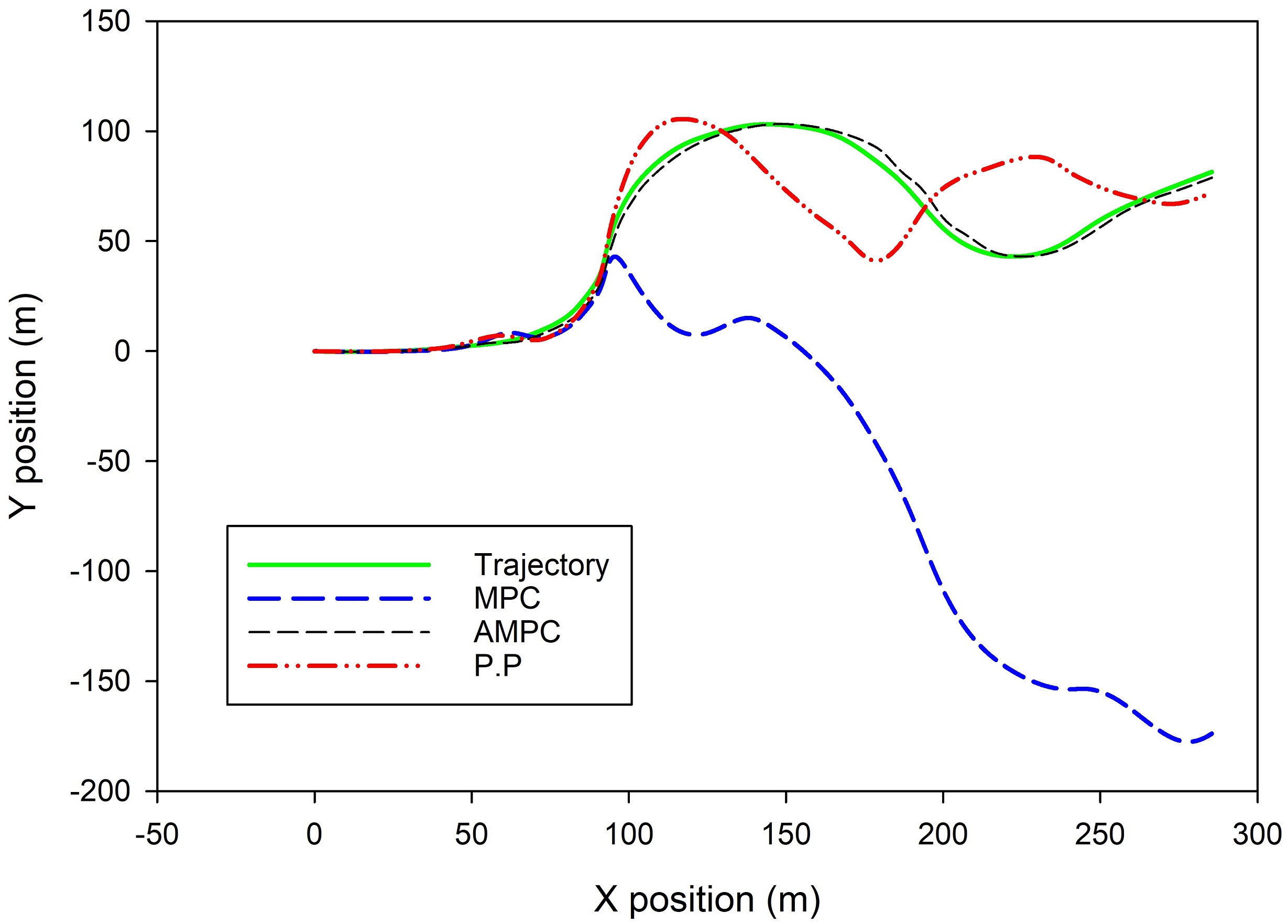}
\caption{Path tracking for \textbf{$Sc_2$} under wind disturbance.}
\label{fig:22}
\end{figure}

\begin{figure}[!h]
\centering
\includegraphics[width=8cm,height=3.9cm]{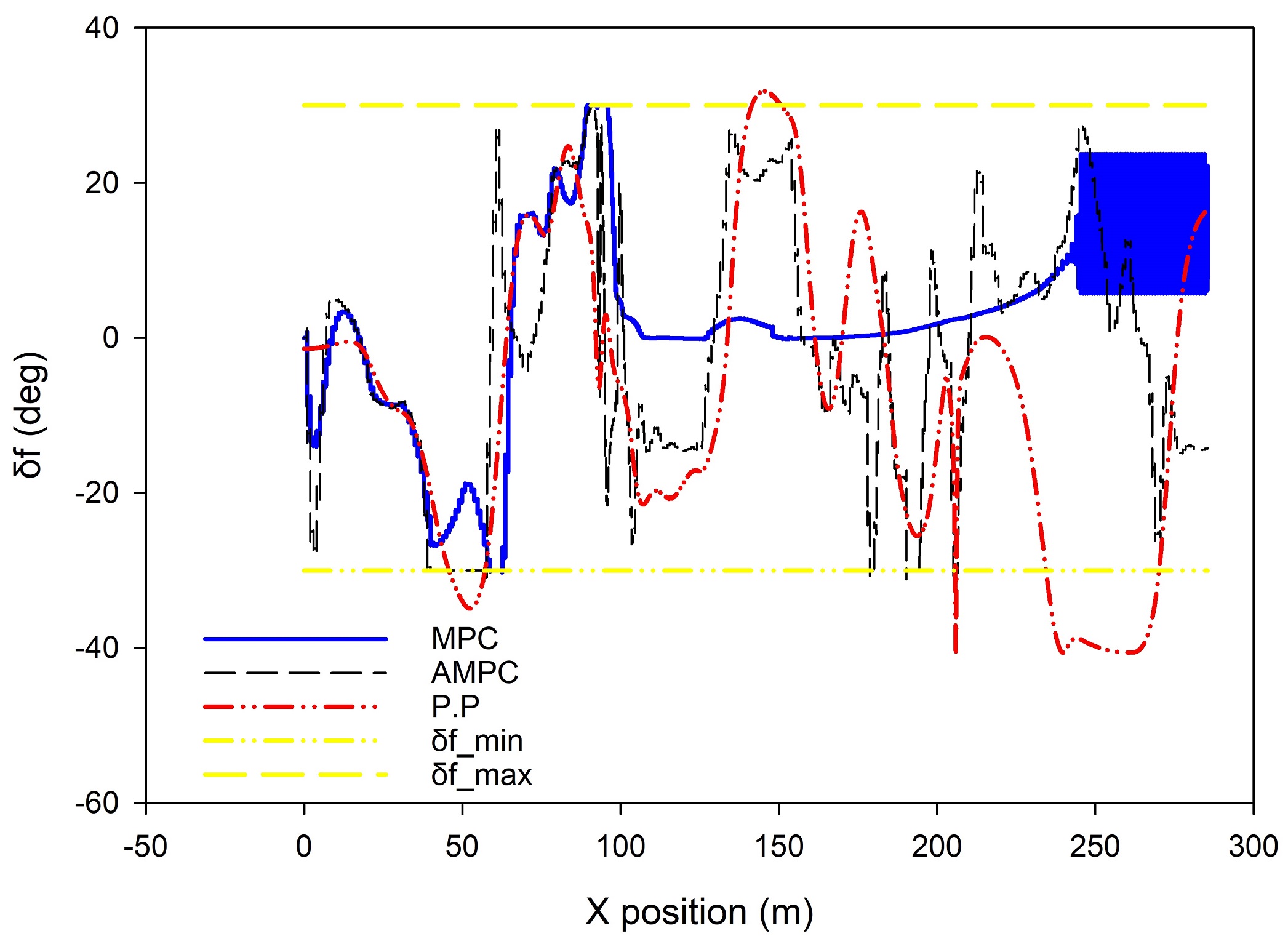}
\caption{Control signal for \textbf{$Sc_3$} under wind disturbance.}
\label{fig:23}
\end{figure}

\begin{figure}[!h]
\centering
\includegraphics[width=8cm,height=3.9cm]{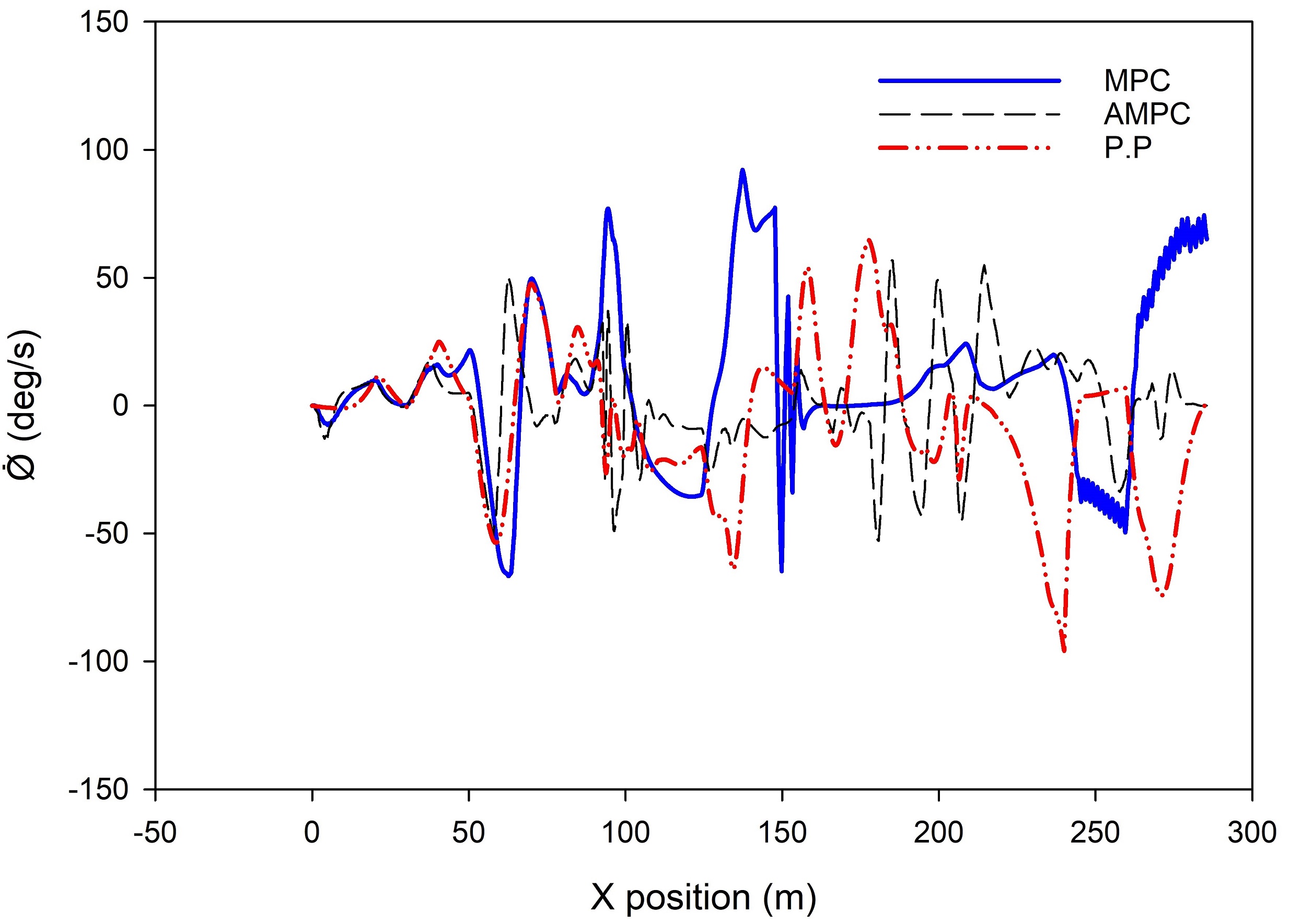}
\caption{Heading rate for \textbf{$Sc_3$} under wind disturbance.}
\label{fig:24}
\end{figure}

For the same trajectory and velocity profile, the lateral wind disturbance profile shown in Fig. \ref{fig:21} is applied to the high fidelity vehicle model to verify the controller robustness. The simulation results (Fig. \ref{fig:22}-\ref{fig:24}) show that at certain wind speeds, the classic MPC is not able to reject the external disturbance and its performance deteriorates and diverges. On the other hand, the proposed AMPC still manages to track the trajectory despite the introduced disturbances. This is due to its adaptive feature with the optimized parameters which allows the AMPC to overcome the perturbations.

\section{Conclusions}

This paper addressed the lateral control task in autonomous vehicles for which an AMPC controller is proposed. A new improved PSO algorithm is also proposed to optimize and tune the AMPC controller for varying working conditions. The MPC controller is formulated with a linearized bicycle model and Laguerre functions and solved through quadratic programming. The parameter optimization is performed for different scenarios, and the generated optimal parameters are used for controller adaptation. The performance of the proposed controller is tested against the classic MPC and the Pure Pursuit controller using a high fidelity vehicle model.  

The results demonstrate that the designed AMPC is able to perform much better tracking accuracy especially for varying longitudinal velocities. Moreover, the AMPC is able to handle varying working conditions and reject external disturbances. These good performances are achieved thanks to its adaptive feature and well optimized parameters through the proposed PSO algorithm. Future research shall address the possibility of using neural networks and adaptive neuro-fuzzy inference systems to learn the optimized parameters by the proposed PSO algorithm for adaptation purposes.


\bibliographystyle{ieeeconf}
\bibliography{library}

\begin{thebibliography}{99}

\bibitem{1} E. Yurtsever, J. Lambert, A. Carballo, and K. Takeda, “A Survey of Autonomous Driving: Common Practices and Emerging Technologies,” IEEE Access, vol. 8, pp. 58443–58469, 2020

\bibitem{2} Q. Yao, Y. Tian, Q. Wang, and S. Wang, “Control Strategies on Path Tracking for Autonomous Vehicle: State of the Art and Future Challenges,” IEEE Access, vol. 8, pp. 161211–161222, 2020

\bibitem{3} S. Bacha, R. Saadi, M. Y. Ayad, A. Aboubou, and M. Bahri, “A review on vehicle modeling and control technics used for autonomous vehicle path following,” Int. Conf. Green Energy Convers. Syst. GECS 2017


\bibitem{4} M. Buehler, K. Iagnemma, and S. Singh, The DARPA urban challenge:
autonomous vehicles in city traffic. Springer, 2009

\bibitem{5} B. Paden, M. Čáp, S. Z. Yong, D. Yershov, and E. Frazzoli, “A survey of motion planning and control techniques for self-driving urban vehicles,” IEEE Trans. Intell. Veh., vol. 1, no. 1, pp. 33–55, 2016

\bibitem{6} M.-W. Park, S.-W. Lee, and W.-Y. Han, ‘Development of lateral control
system for autonomous vehicle based on adaptive pure pursuit algorithm,’ in Proc. 14th Int. Conf. Control, Autom. Syst. (ICCAS), Oct. 2014, pp. 1443–1447.

\bibitem{7} J.X. Wang, M.M. Dai, G.D. Yin, and N. Chen Output-feedback robust control for vehicle path tracking considering different human drivers’ characteristics. Mechatronics 2018; 50: 402–412

\bibitem{8} G. Han, W. Fu, W. Wang, and Z. Wu, “The lateral tracking control for the intelligent vehicle based on adaptive PID neural network,” Sensors (Switzerland), vol. 17, no. 6, pp. 1–15, 2017

\bibitem{26} Y. Kebbati, N. Ait-oufroukh, V. Vigneron, D. Ichalal, and D. Gruyer, "Optimized self-adaptive PID speed control for autonomous vehicles," ICAC 2021 - 2021 26th IEEE Int. Conf. Autom. Comput. no. September, pp. 81–86, 2021.

\bibitem{9} C. Hu, R. Wang, and F. Yan, “Integral Sliding Mode-Based Composite Nonlinear Feedback Control for Path Following of Four-Wheel Independently Actuated Autonomous Vehicles,” IEEE Trans. Transp. Electrif., vol. 2, no. 2, pp. 221–230, 2016


\bibitem{10} B. Zhang, C. Zong, G. Chen, and B. Zhang, “Electrical Vehicle Path Tracking Based Model Predictive Control with a Laguerre Function and Exponential Weight,” IEEE Access, vol. 7, pp. 17082–17097, 2019

\bibitem{11} H. Guo, D. Cao, H. Chen, Z. Sun, and Y. Hu, “Model predictive path following control for autonomous cars considering a measurable disturbance: Implementation, testing, and verification,” Mech. Syst. Signal Process., vol. 118, pp. 41–60, 2019

\bibitem{12} H. Wang, B. Liu, X. Ping, and Q. An, “Path Tracking Control for Autonomous Vehicles Based on an Improved MPC,” IEEE Access, vol. 7, pp. 161064–161073, 2019

\bibitem{13} M. S. Akbari, A. A. Safavi, N. Vafamand, T. Dragicevic, and J. Rodriguez, “Fuzzy Mamdani-based Model Predictive Load Frequency Control,” 2020 IEEE 11th Int. Symp. Power Electron. Distrib. Gener. Syst. PEDG 2020, pp. 7–12, 2020


\bibitem{14} V. Ramasamy, R. K. Sidharthan, R. Kannan, and G. Muralidharan, “Optimal tuning of model predictive controller weights using genetic algorithm with interactive decision tree for industrial cement kiln process,” Processes, vol. 7, no. 12, 2019


\bibitem{15} H. Moumouh, N. Langlois, and M. Haddad, “A Novel Tuning approach for MPC parameters based on Artificial Neural Network,” IEEE Int. Conf. Control Autom. ICCA, vol. 2019-July, pp. 1638–1643, 2019

\bibitem{16} F. Lin, Y. Chen, Y. Zhao, and S. Wang, “Path tracking of autonomous vehicle based on adaptive model predictive control,” Int. J. Adv. Robot. Syst., vol. 16, no. 5, pp. 1–12, 2019

\bibitem{17} M. Bujarbaruah, X. Zhang, H. E. Tseng, and F. Borrelli, “Adaptive MPC for Autonomous Lane Keeping,” 2018

\bibitem{18} L. Wang,Model  Predictive  Control  System  Design  and  Implementation Using MATLAB, vol. 53.  2016.

\bibitem{19} B. Song, Z. Wang, and L. Zou, “An improved PSO algorithm for smooth path planning of mobile robots using continuous high-degree Bezier curve,” Appl. Soft Comput., vol. 100, p. 106960, 2021


\bibitem{20}  Y. Xie and W. Zhang, “A novel tuning method for PID controller based on improved PSO algorithm for unstable plants with time delay,” Chinese Control Conf. CCC, vol. 2019-July, pp. 4270–4275, 2019

\bibitem{21} Z. Tian, Y. Ren, and G. Wang, “Short-term wind speed prediction based on improved PSO algorithm optimized EM-ELM,” Energy Sources, Part A Recover. Util. Environ. Eff., vol. 41, no. 1, pp. 26–46, 2019

\bibitem{22} Y. Shi and R. Eberhart, “Modified particle swarm optimizer,” Proc. IEEE Conf. Evol. Comput. ICEC, no. February 2015, pp. 69–73, 1998


\bibitem{23} Z. H. Zhan, J. Zhang, Y. Li, and H. S. H. Chung, “Adaptive particle swarm optimization,” IEEE Trans. Syst. Man, Cybern. Part B Cybern., vol. 39, no. 6, pp. 1362–1381, 2009


\bibitem{24} Y. Del Valle, G. K. Venayagamoorthy, S. Mohagheghi, J. C. Hernandez, and R. G. Harley, “Particle swarm optimization: Basic concepts, variants and applications in power systems,” IEEE Trans. Evol. Comput., vol. 12, no. 2, pp. 171–195, Apr. 2008

\bibitem{25} H.B. Pacejka,  (2008). Vehicle System Dynamics : International Journal of Vehicle Mechanics and Mobility. International Journal of Vehicle Mechanics and Mobility, (August 2012), 37–41

\end{thebibliography}

\end{document}